\documentclass[12pt,a4paper]{article}

\usepackage[british]{babel}

\usepackage[T1]{fontenc}

\usepackage[utf8]{inputenc}

\usepackage{csquotes}

\usepackage[backend=biber, maxbibnames=10, maxcitenames=4, maxalphanames=4, bibencoding=utf8, style=alphabetic, isbn=false, url=false, doi=true]{biblatex}
\renewbibmacro{in:}{}
\renewbibmacro*{doi+eprint+url}{%
	\printfield{doi}%
	\newunit\newblock%
	\iftoggle{bbx:eprint}{%
		\usebibmacro{eprint}%
	}{}%
	\newunit\newblock%
	\iffieldundef{doi}{%
		\usebibmacro{url+urldate}}%
	{}%
}

\usepackage{setspace}

\usepackage{lmodern}

\usepackage[indentafter]{titlesec}
\titleformat{name=\section}{}{\thetitle.}{0.8em}{\centering\scshape}
\titleformat{name=\subsection}[runin]{}{\thetitle.}{0.5em}{\bfseries}[.]
\titleformat{name=\subsubsection}[runin]{}{\thetitle.}{0.5em}{\itshape}[.]
\titleformat{name=\paragraph,numberless}[runin]{}{}{0em}{}[.]
\titlespacing{\paragraph}{0em}{0em}{0.5em}
\titleformat{name=\subparagraph,numberless}[runin]{}{}{0em}{}[.]
\titlespacing{\subparagraph}{0em}{0em}{0.5em}

\usepackage{url}

\usepackage{mathtools}

\usepackage{amssymb}

\usepackage{amsthm}

\usepackage{mathrsfs}

\usepackage{dsfont}

\usepackage{enumitem}
\setlist[enumerate]{noitemsep, partopsep=0pt, topsep=0pt, parsep=0pt, itemsep=0pt}
\setlist[itemize]{noitemsep, partopsep=0pt, topsep=0pt, parsep=0pt, itemsep=0pt}

\usepackage{interval}

\intervalconfig{soft open fences}

\usepackage[normalem]{ulem}

\usepackage[textsize=tiny]{todonotes}

\setuptodonotes{fancyline}
\usepackage[stretch=10]{microtype}

\usepackage[affil-it]{authblk}

\usepackage{xfrac}

\graphicspath{./img}

\usepackage{hyperref}

\usepackage[nameinlink,noabbrev,capitalize]{cleveref}

\usepackage[plain]{fullpage}

\newlist{theoenum}{enumerate}{1} 
\setlist[theoenum]{label=\normalfont(\roman*), ref=\theproposition~\normalfont(\roman*), noitemsep, partopsep=0pt, topsep=0pt, parsep=0pt, itemsep=0pt}
\crefalias{theoenumi}{theorem}

\newlist{maintheoenum}{enumerate}{1} 
\setlist[maintheoenum]{label=\normalfont(R\arabic*), ref=\theproposition~\normalfont(R\arabic*), noitemsep, partopsep=0pt, topsep=0pt, parsep=0pt, itemsep=0pt}
\crefalias{maintheoenumi}{~} 

\theoremstyle{plain}
\newtheorem{theorem}{Theorem}
\newtheorem*{theorem*}{Theorem}
\newtheorem{proposition}[theorem]{Proposition}
\newtheorem{corollary}[theorem]{Corollary}

\newtheorem{lemma}[theorem]{Lemma}
\newtheorem*{lemma*}{Lemma}

\theoremstyle{definition}
\newtheorem{definition}[theorem]{Definition}

\theoremstyle{remark}

\newtheorem{remark}[theorem]{Remark}
\newtheorem*{remark*}{Remark}			

\newcommand*\diff{\mathop{}\!\mathrm{d}}
\newcommand{\R}{\mathds{R}}
\newcommand{\Z}{\mathds{Z}}

\newcommand{\He}{\mathds{H}}
\newcommand{\E}{\mathrm{E}}
\newcommand{\J}{\mathrm{J}}
\newcommand{\D}{\mathcal{D}} 
\newcommand{\A}{\mathscr{A}}

\newcommand{\me}{\mathrm{e}}

\renewcommand{\L}{\mathrm{L}}
\renewcommand{\exp}{\mathrm{exp}}
\renewcommand{\epsilon}{\varepsilon}

\DeclareMathOperator{\Lie}{\mathrm{Lie}}
\DeclareMathOperator{\T}{\mathrm{T}}
\DeclareMathOperator{\vspan}{\mathrm{span}}
\DeclareMathOperator{\Ima}{\mathrm{Im}}
\DeclareMathOperator{\Kern}{\mathrm{Ker}}

\date{\footnotetext[0]{Date: \today.}}
\title{\textbf{\uppercase{\large{Regularity and continuity properties\\   of the sub-Riemannian exponential map}}}}
\author{Samuël Borza\thanks{\href{mailto:sborza@sissa.it}{sborza@sissa.it}} }
\affil[]{International School for Advanced Studies (SISSA, Trieste)}
\author{Wilhelm Klingenberg\thanks{\href{mailto: wilhelm.klingenberg@dur.ac.uk}{wilhelm.klingenberg@dur.ac.uk }} }
\affil[]{Department of Mathematics, Durham University}					

\begin{document}

\maketitle							

\providecommand{\keywords}[1]
{
	\textbf{\textit{Keywords---}} #1
}

\providecommand{\msc}[1]
{
	\textbf{\textit{MSC (2020)---}} #1
}

\vspace{-1.5cm}\begin{abstract}
	 We prove a version of Warner's regularity and continuity properties for the sub-Riemannian exponential map. The regularity property is established by considering sub-Riemannian Jacobi fields while the continuity property follows from studying the Maslov index of Jacobi curves. We finally show how this implies that the exponential map of the three dimensional Heisenberg group is not injective in any neighbourhood of a conjugate vector.
\end{abstract}
\keywords{Conjugate points, Sub-Riemannian geometry, Metric geometry}\\
\msc{53C17, 53B15, 53B99, 49J15}				

\tableofcontents

\section{Introduction} \label{intro}

In his work \cite{warner1965} on the conjugate locus in Riemannian geometry, Warner introduced the notion of \textit{regular exponential map}, a map $F : \T_q(M) \to M$, where $M$ is a finite dimensional smooth manifold, that satisfies three conditions: a non-vanishing speed condition along rays, a regularity and a continuity condition. Furthermore, Warner shows that such a map is non-injective in any neighbourhood of any singularities of $F$. This is done through studying the normal forms of $F$ around singularities, namely the points where the jacobian determinant of $F$ vanish. Warner then proves that the exponential map of a Finsler manifold is regular in this sense, giving an alternative proof of a result due to Morse and Littauer \cite{morse}. Warner's theorem was adapted to Lorentzian structures in \cite{rosquist1983} and then to semi-Riemannian manifolds in \cite{szeghy2008}.

In the present work, we adapt Warner's conditions for the exponential map in sub-Riemannian geometry. Because of the lack of a Levi-Civita connection, the study of geodesics is carried out from the \textit{Hamiltonian viewpoint} (see \cref{prelim} for a summary of the theory). Length minimisers are found to be \textit{normal} and/or \textit{abnormal}. Normal geodesics are solutions of a Hamiltonian system of differential equations with initial conditions taking values in the cotangent space. The sub-Riemannian exponential map is the projection of the corresponding Hamiltonian flow. \textit{Abnormal} geodesics can also appear: they are length minimisers that satisfy a condition not characterised by a differential equation. Strongly normal geodesics are those that do not contain abnormal subsegments.

Let us introduce some notations to state our main results. The normal geodesic starting at $q \in M$ with initial covector $\lambda_0 \in \mathrm{T}^*_q(M)$ is denoted by $\gamma_{q, \lambda_0}$ and $I_{q, \lambda_0}$ is its maximal domain. The ray through $\lambda_0$ is the map defined by $r_{q, \lambda_0}(t) := t \lambda_0$ for $t \in I_{q, \lambda_0}$. We denote by $H_q$ the restriction of the Hamiltonian to a fiber $\mathrm{T}^*_q(M)$. For $A \in  \Kern(\diff_{\lambda_0} \exp_q)$, the sub-Riemannian Jacobi field $J_A$ along $\gamma_{q, \lambda_0}$ is the one with initial values $(0, A)$. Choosing a symplectic moving frame along $\gamma_{q, \lambda_0}$ allows us to introduce $\nabla J_A$, a (non canonical) derivative of $J_A$ along the curve. The theory of sub-Riemannian Jacobi fields will be detailled in \cref{jacobifields}.

\begin{theorem}[Regularity of the sub-Riemannian exponential map]
	\label{regsubriemexp}
	Let $M$ be a sub-Riemannian manifold, $q \in M$ and $\exp_q : \mathscr{A}_q \to M$ be the corresponding exponential map with domain $\A_q \subseteq \T^*_q(M)$. Then,
	\begin{maintheoenum}
		\item \label{(R1)} The map $\mathrm{exp}_q$ is $\mathcal{C}^\infty$ on $\A_q$ and, for all $\lambda_0 \in \A_q \setminus H^{-1}_q(0)$ and all $t \in I_{q, \lambda_0}$, we have $\mathrm{d}_{t \lambda_0} \mathrm{exp}_q(\mathrm{\dot{r}}_{q, \lambda_0}(t)) \neq 0_{\mathrm{exp}_q(t \lambda_0)}$;
		\item \label{(R2)} For every $\lambda_0 \in \A_q \setminus \{ 0 \}$ and every symplectic moving frame along the cotangent lift $\lambda(t)$ of the normal geodesic $\gamma(t):= \exp_q(t \lambda_0)$, the map 
		\[
		\Kern(\diff_{\lambda_0} \exp_q) \to \T_{\exp_q(\lambda_0)}(M)/\diff_{\lambda_0} \exp_q(\T_{\lambda_0}(\T^*_q(M))),
		\]
		sending $A$ to $\nabla J_A(1) + \diff_{\lambda_0} \exp_q(\T_v(\T^*_q(M)))$, is a linear isomorphism;
		\item \label{(R3)} Let $\lambda_0 \in \mathscr{A}_q \setminus H_q^{-1}(0)$ be a covector such that the corresponding geodesic $\gamma(t) := \exp_q(t \lambda_0)$ is strongly normal. Then, there exists a radially convex neighbourhood $\mathcal{V}$ of $\lambda_0$ such that for every ray $r_{q, \overline{\lambda}_0}$ which intersects $\mathcal{V}$ that does not contain abnormal subsegments in $\mathcal{V}$, the number of  singularities of $\exp_q$ (counted with multiplicities) on $\mathrm{Im}(r_{q, \overline{\lambda}_0}) \cap \mathcal{V}$ is constant and equals the order of $\lambda_0$ as a singularity of $\exp_q$, i.e. $\mathrm{dim} (\mathrm{Ker}(\diff_{\lambda_0} \exp_{q}))$.
	\end{maintheoenum}
\end{theorem} 
Condition (R1) follows from the constant speed property of normal geodesics, see \cref{normalextremals}. The rich theory of sub-Riemannian Jacobi fields, developed in \cite{subriemjacobi} for example, will help us to prove condition (R2) in \cref{jacobifields}. In Riemannian geometry, condition (R3) is a consequence of Morse's index theory. His ideas are adapted to this context with the Maslov index and the condition (R3) will be obtained in \cref{R3}. 

Warner uses these three conditions in \cite{warner1965} to conclude that the Riemannian exponential map is not locally injective around a singularity. This result is originally due to Morse and Littauer \cite{morse}. In this paper, Warner's method is used for the first time in sub-Riemannian geometry, in the case of the three dimensional Heisenberg group $\mathds{H}$ in \cref{examples}. We expect it to work for larger classes of sub-Riemannian manifolds as well.

\begin{theorem}
	\label{thm:noninjH}
	For $q \in \mathds{H}$, the sub-Riemannian exponential map $\exp_q : \mathscr{A}_q \to \mathds{H}$ is not injective on any neighbourhood of a conjugate vector $\lambda_0 \in \mathscr{A}_q \setminus H_q^{-1}(0)$.
\end{theorem}


\section*{Acknowledgements}

This work was supported by the UK Engineering and Physical Sciences Research Council (EPSRC) grant EP/N509462/1 (1888382) for Durham University. 

This project has received funding from the European Research Council (ERC) under the European Union’s Horizon 2020 research and innovation programme (grant agreement No. 945655).		

\section{Preliminaries} \label{prelim}

In this section, we set up the basics of sub-Riemannian geometry. We rely on \cite{comprehensive2020} and \cite{agrachev2008} for the general theory.

\subsection{Sub-Riemannian geometry}

We begin with a definition of \emph{sub-Riemannian structure}.

\begin{definition}
	A triple $(E, \langle \cdot, \cdot \rangle_E, f_E)$ induces a sub-Riemannian structure on a manifold $M$ if 
	\begin{enumerate}[label=\normalfont(\roman*)]
		\item $E$ is a vector bundle over $M$,
		\item $\langle \cdot, \cdot \rangle_E$ is an inner product on $E$,
		\item $f_E : E \to \T(M)$ is a morphism of vector bundles.
	\end{enumerate}
\end{definition}
The family $\D$ of \textit{horizontal vector fields} is defined as $\D := \{ f_E \circ u \mid u \text{ is a section of } E \}$ and the distribution at a point $q \in M$ is $\D_q := \{ v(x) \mid v \in \D \}$. The \textit{rank} of the sub-Riemannian structure at $q \in M$ is $\mathrm{rank}(q) := \mathrm{dim}(\D_q)$. Observe that in our definition, a sub-Riemannian manifold can be rank-varying.

\begin{definition}
	We say that curve $\gamma : [0, T] \to M$ is \textit{horizontal} if $\gamma$ is Lipschitz in charts and if there exists $u \in \mathrm{L}^2(\interval{0}{T}, E)$ such that for almost every $t \in [0, T]$, we have $u(t) \in E_{\gamma(t)}$ and $\dot{\gamma} (t) = f_E(u(t))$.
	The \textit{sub-Riemannian length} and the \textit{sub-Riemannian energy} of $\gamma$ are defined by
	\begin{equation}
		\label{length}
		\L(\gamma) = \int_{0}^{T} \lVert \dot{\gamma} (t) \rVert_{\D_{\gamma(t)}} \diff t, \ \J(\gamma) = \dfrac{1}{2} \int_{0}^{T} \lVert \dot{\gamma} (t) \rVert^2_{\D_{\gamma(t)}} \diff t
	\end{equation}
	where $\lVert v \rVert_{\D_q} := \min \left\{ \sqrt{\langle u, u \rangle_{E_q}} \mid u \in E_q \text{ and } f_E(u) = (q, v) \right\}$ for $v \in \D_q$ and $q \in M$.
\end{definition}

\begin{remark}
	The norm $\lVert \cdot \rVert_{\D_q}$ is well defined, induced by an inner product $\langle \cdot, \cdot \rangle_{\D_q}$ via the polarisation formula and the map $t \mapsto \lVert \dot{\gamma} (t) \rVert_{\D_{\gamma(t)}}$ is measurable.
\end{remark}

In the case where every two points can be joined by a horizontal curve, we have a well defined distance function on $M$.

\begin{definition}
	The distance of a sub-Riemannian manifold $M$, also called the \textit{Carnot-Carathéodory distance}, is defined by
	\[
	\diff(x, y) := \inf \{ \L(\gamma) \mid \gamma : [0,T] \to M \text{ is horizontal and } \gamma(0) = x \text{ and } \gamma(T) = y \}.
	\]
\end{definition}

In this work, we assume that the sub-Riemannian structures satisfy the \textit{Hörmander condition}, that is to say $\Lie_q(\D) = \T_q(M)$ for all $q \in M$. We also say that $\D$ is \textit{bracket-generating}. This is motivated by the following well-known result. 

\begin{theorem}[Chow–Rashevskii theorem, see {\cite[Theorem 3.31.]{comprehensive2020}}]
	Let $M$ be a sub-Riemannian manifold such that its distribution $\D$ is $\mathcal{C}^\infty$ and satisfies the Hörmander condition. Then, $(M, \diff)$ is a metric space and the manifold and metric topology of $M$ coincide.
\end{theorem}


Given $m$ global vector fields $X_1, \dots, X_m : M \to \T(M)$ on a  manifold $M$, we can build on it a sub-Riemannian structure in the following way. We set $E := M \times \mathds{R}^m$ as the trivial bundle of rank $m$, $f_E : E \to \T(M) : (q, (u_1, \dots, u_m)) \mapsto \sum_{k = 1}^{m} u_k X_k(q)$ and finally the metric on $E$ is the Euclidean one. In this way, we induce an inner product on $\D_q = \vspan\{X_1(q), ..., X_m(q)\}$ by the polarization formula applied to the norm
\[
\|u\|^2_{\D_q} := \min \left\{ \sum_{k = 1}^{m} u^2_i \ \middle| \ \sum_{k = 1}^{m} u_i X_k(q) = u \right\}.
\]
The family $(X_1, \dots, X_m)$ is said to be a \textit{generating family} of the sub-Riemannian manifold. A \textit{free} sub-Riemannian structure is one that is induced from a generating family. Every sub-Riemannian structure is \textit{equivalent} to a free one (see \cite[Section 3.1.4]{comprehensive2020}). From now on, we will therefore assume, without loss of generality, that every sub-Riemannian manifold is free.

\subsection{End-point map and Lagrange multipliers}

We fix a sub-Riemannian manifold $M$ for which the family $(X_1, \dots, X_m)$ is generating. From an optimal control point of view, a curve $\gamma : \interval{0}{T} \to M$ with initial value $\gamma(0) = q \in M$ is horizontal if there exists $u \in \L^{2}(\interval{0}{T}, \R^m)$, called a \textit{control}, such that $\dot{\gamma}(t) = \sum_{k = 0}^m u_k(t) X_k(\gamma(t))$ for almost every $t \in \interval{0}{T}$.

In fact, by Carathéodory's theorem for ordinary differential equations, we know that there exists a unique maximal Lipschitz solution to the Cauchy problem
\begin{equation}
\label{cauchyprob}
\begin{cases}
	\dot{\gamma}(t) = \sum_{k = 0}^m u_k(t) X_k(\gamma(t)) \\
	\gamma(t_0) = q
\end{cases}
\end{equation}
for every $u \in \L^{2}(\interval{0}{T}, \R^m)$, $q \in M$ and $t_0 \in \interval{0}{T}$. We denote such a solution by $\gamma_{t_0, q, u}$ and we can now introduce the end-point map. 

\begin{definition}
	Let $q \in M$ and $T > 0$. The \textit{end-point map at time} $T > 0$ is the smooth map
	\[
	\E_{q, T} : \mathcal{U} \to M : u \mapsto \gamma_{0, q, u}(T),
	\]
	where $\mathcal{U} \subseteq \L^{2}(\interval{0}{T}, \R^m)$ is the open subset of controls such that $\gamma_{0, q, u}$, solution to the Cauchy problem \eqref{cauchyprob}, is defined on the whole interval $\interval{0}{T}$.
\end{definition}


On the space of controls $\L^{2}(\interval{0}{T}, \R^m)$, we can define a \textit{length functional}, as well as a corresponding \textit{energy functional}
\[
L(u) := \int_{0}^{T} \| u(t) \|_{\R^m} \diff t,  \ J(u) := \dfrac{1}{2} \int_{0}^{T} \| u(t) \|^2_{\R^m} \diff t.
\]
Given an horizontal curve $\gamma : \interval{0}{T} \to M$, we define at every differentiability point of $\gamma$ the \textit{minimal control} $\overline{u}$ associated with $\gamma$
\[
\overline{u}(t) := \arg \min \left\{ \| u \|_{\R^m} \mid u \in \R^m, \dot{\gamma}(t) = \sum_{k = 0}^m u_k X_k(\gamma(t)) \right\}.
\]
In particular, the previous functionals are related to the sub-Riemannian length and energy: $\L(\gamma) = L(\overline{u})$ and $\J(\gamma) = J(\overline{u})$, where $\overline{u}$ is the minimal control associated with $\gamma$. Through the Cauchy problem \eqref{cauchyprob}, it is clear that finding a length minimiser for $\L$ among the horizontal curves with fixed end-points $\gamma(0) = q$ and $\gamma(T) = q'$ is equivalent to finding a minimal control for $L$ for which the associated path joins $q$ and $q'$. Furthermore, we have the following classical correspondence: a horizontal curve $\gamma : \interval{0}{T} \to M$ joining $q$ to $q'$ is a minimiser of $\J$ if and only if it is a minimiser of $\L$ and is parametrised by constant speed.

In terms of the end-point map, the problem of finding the minimisers joining two fixed point $q, q' \in M$ is thus equivalent to solving the constrained variational problem
\begin{equation}
	\label{constraint}
	\min \left\{J(u) \ \middle| \ u \in \E^{-1}_{q, T}(q') \right\}.
\end{equation}
The Lagrange multipliers rule provides a necessary condition to be satisfied by a control $u$ that is a constrained critical point for \eqref{constraint}.
\begin{proposition}[Lagrange Multipliers]
	Let $u \in \mathcal{U}$ be an optimal control for the variation problem \eqref{constraint}. Then at least one of the following statements holds true
	\begin{enumerate}[label=\normalfont(\roman*)]
		\item \label{lagrangenormal} there exists $\lambda(T) \in \T^*_q(M)$ such that $\lambda(T) \circ D_u \E_{q, T} = \diff_u J$;
		\item \label{lagrangeabnormal} there exists $\lambda(T) \in \T^*_q(M) \setminus \{0\}$ such that $\lambda(T) \circ D_u \E_{q, T} = 0$.
	\end{enumerate}	
\end{proposition}

An optimal control is called \textit{normal} (resp. \textit{abnormal}) when it satisfies the condition \ref{lagrangenormal} (resp. \ref{lagrangeabnormal}). The same terminology is used for the corresponding curve $\gamma_u$. We note that the extremal trajectory $\gamma_u$ could be both normal and abnormal. A normal trajectory $\gamma : \interval{0}{T} \to M$ is called \textit{strictly normal} if it is not abnormal. If, in addition, the restriction $\gamma|_{\interval{0}{s}}$ is strictly normal for every $s > 0$, we say that $\gamma$ is \textit{strongly normal}. It can be seen that $\gamma$ is strongly normal if and only if the normal geodesic $\gamma$ does not contain any abnormal segment. 

\subsection{Characterisation of sub-Riemannian geodesics}

Now that we can turn a sub-Riemannian manifold into a metric space, we would like to study the geodesics associated with its distance function. These would be horizontal curves that are locally minimising the sub-Riemannian length functional. Because of the lack of a torsion-free metric connection, we can not have a geodesic equation through a covariant derivative. Rather, sub-Riemannian geodesic  are characterised via Hamilton's equations.

We recall that the Hamiltonian vector field of a map $a \in C^{\infty}(\T^*(M))$ is the unique vector field $\overrightarrow{a}$ on $\T^*(M)$ that satisfies
\[
\omega(\cdot, \overrightarrow{a}(\lambda)) = \diff_\lambda a, \qquad \forall \lambda \in \T^*(M),
\]
where $\omega$ denotes the canonical symplectic form on the cotangent bundle $\T^*(M)$. \textbf{}

The smooth \textit{control-dependent Hamiltonian} of a sub-Riemannian structure is the map $h : \R^m \times \T^*(M) \to \R$ defined as 
\[
h_u(\lambda) = \sum_{k = 1}^{m} u_k \langle \lambda, X_k(\pi(\lambda)) \rangle - \dfrac{1}{2} \sum_{k = 1}^{m} u_k ^2.
\]
It is easy to see that, by strict convexity, there exists a unique maximum $\overline{u}(\lambda)$ of $u \mapsto h_u(\lambda)$ for every $\lambda \in \T^*(M)$. Therefore, a \textit{maximized Hamiltonian}, or simply \textit{Hamiltonian}, is well defined
\[
H : \T^*(M) \to \mathds{R} : \lambda \mapsto H(\lambda) := \max_{u \in \R^m} h_u(\lambda).
\]
Furthermore, the Hamiltonian $H$ can written in terms of the generating family of the sub-Riemannian structure $(X_1, \dots, X_m)$, as follows
\[
H(\lambda) = \frac{1}{2}\sum_{k = 1}^m h_k(q, \lambda_0)^2,
\]
where $\lambda := (q, \lambda_0) \in T^*(M)$ and $h_k(q, \lambda_0) := \langle \lambda_0, X_k(q) \rangle$. For $q \in M$, we will also write $H_q$ for the restriction of $H$ to the cotangent space $\T_q^*(M)$.

The Lagrange multipliers rule may be further developed to characterise normal extremals as curves that satisfy Hamilton's differential equation. Alternatively, the following result can also be seen as an application of the Pontryagin Maximum Principle to the sub-Riemannian length minimisation problem.

\begin{theorem}[Pontryagin's Maximum Principle]
	Let $\gamma : \interval{0}{T} \to M$ be a horizontal curve which is a length minimiser among horizontal curves, and parametrised by constant speed. Then, there exists a Lipschitz curve $\lambda(t) \in \mathrm{T}^*_{\gamma(t)}(M)$ such that one and only one of the following is satisfied:
	\begin{enumerate}
		\item[(N)] \label{hamiltoneq} $\dot{\lambda} = \overrightarrow{H}(\lambda)$, where $\overrightarrow{H}$ is the unique vector field in $\mathrm{T}^*(M)$ such that $\omega(\cdot, \overrightarrow{H}(\lambda)) = \mathrm{d}_\lambda H$ for all $\lambda \in \mathrm{T}^*(M)$;
		\item[(A)] \label{abnormalcondi} $\omega_{\lambda(t)}(\dot{\lambda} (t), \cap_{k = 1}^m \mathrm{ker}(\mathrm{d}_{\lambda(t)} h_k)) = 0$ for all $t \in \interval{0}{T}$.
	\end{enumerate}
	Moreover, in the case (A), we have $\lambda(t) \neq 0$ for every $t \in \interval{0}{T}$.
\end{theorem}

If $\lambda$ satisfies $(N)$ (resp. $(A)$), we will also say that $\lambda$ is a normal extremal (resp. abnormal extremal). The projection of a normal extremal to $M$ is locally minimising, that is to say it is a normal geodesic parametrised by constant-speed. This can be seen of an application of the invariance of the Hilbert integral around small subsegments of the trajectory, defined through the Poincarré-Cartan one-form (see \cite[Section 4.7]{comprehensive2020}). However, the projection of an abnormal extremal to $M$ might not be locally minimising. 

It does happen that a sub-Riemannian structure does not have any non-trivial (i.e. non-constant) abnormal geodesic (the trivial geodesic is always abnormal as soon as $\mathrm{rank}(\D_q) < \mathrm{dim} (M)$). In this case, the sub-Riemannian manifold is said to be \emph{ideal}.


The theory of ordinary differential equations provides the existence of a maximal solution to (N) for every given an initial condition $(q, \lambda_0) \in \T^*(M)$. The flow of Hamilton's equation is denoted by $\me^{t \overrightarrow{H}}$. We finally turn our attention to the central object of this work: the sub-Riemannian \textit{exponential map}.

\begin{definition}
	The sub-Riemannian \textit{exponential map} at $q \in M$ is the map
	\[
	\exp_q : \A_q \to M : \lambda \mapsto \pi(\me^{\overrightarrow{H}}(\lambda)),
	\]
	where $\A_q \subseteq \T^*_q(M)$ is the open set of covectors such that the corresponding solution of (N) is defined on the whole interval $\interval{0}{1}$.
\end{definition}

The sub-Riemannian exponential map $\exp_q$ is smooth. If $\lambda : \interval{0}{T} \to \T^*(M)$ is the normal extremal that satisfies the initial condition $(q, \lambda_0) \in \T^*(M)$, then the corresponding normal extremal path $\gamma(t) = \pi(\lambda(t))$ by definition satisfies $\gamma(t) = \exp_q(t \lambda_0)$ for all $t \in \interval{0}{T}$. If $M$ is complete for the Carathéodory distance, then $\A_q = \T^*_q(M)$ and if in addition there are no stricly abnormal length minimisers, the exponential map $\exp_q$ is surjective. Contrary to the Riemannian case, the sub-Riemannian exponential map is not necessarily a diffeomorphism of a small ball in $\T^*_q(M)$ onto a small geodesic ball in $M$. In fact, $\Ima (\diff_0 \exp_q) = \D_q$ and $\exp_q$ is a local diffeomorphism at $0 \in \T^*_q(M)$ if and only if $\D_q = \T^*_q(M)$. Our aim is to prove that $\exp_q$ is \textit{regular} in the sense of Warner \cite{warner1965}.			

\section{Regularity and Continuity of the sub-Riemannian exponential map} \label{R1}

\subsection{Normal extremals} \label{normalextremals}

As pointed out in the previous section, the normal geodesic $\gamma(t)$ of the sub-Riemannian manifold $M$, with initial point $\gamma(0) = p$ and initial covector $\lambda_0 \in \A_q$ is the projection of the normal extremal $\lambda : I_{q, \lambda_0} \to \T^*(M)$, the solution to Hamilton's geodesic equation with initial value $(q, \lambda_0)$.

The \textit{ray} in $\A_q$ through $\lambda_0$ is the map
\[
r_{q, \lambda_0} : I_{q, \lambda_0} \to \mathrm{T}^*_q(M) : t \mapsto t \lambda_0
\]
where $I_{q, \lambda_0} \subseteq \R^+$ is the maximal interval containing 0 such that $t \lambda_0 \in \A_q$ for every $t \in I_{q, \lambda_0}$. In this way, $\mathrm{\dot{r}}_{q, \lambda_0}(t) \in \T_{t \lambda_0}(\T^*_q(M))$ and identifying $\T_{t \lambda_0}(\T^*_q(M))$ with $\T^*_q(M)$ in the usual way, we have $\mathrm{\dot{r}}_{q, \lambda_0}(t) = \lambda_0$ for every $t \in I_{q, \lambda_0}$.

\begin{proposition}[see \cite{comprehensive2020}, Theorem 4.25]
	Let $\lambda : \interval{0}{T} \to \T^*(M)$ be a normal extremal, that is a solution to Hamilton's equation $\dot{\lambda} = \overrightarrow{H}(\lambda)$. The corresponding normal geodesic $\gamma(t) = \pi(\lambda(t))$ has constant speed and $\frac{1}{2} \| \dot{\gamma} (t) \|^2 = H(\lambda(t))$ for every $t \in \interval{0}{T}$.
\end{proposition} 

\begin{proof}
	The Hamiltonian is constant along a normal trajectory:
	\[
	\dfrac{\diff}{\diff t} H(\lambda(t)) = \diff_\lambda H \circ \dot{\lambda}(t) = \omega(\dot{\lambda}(t), \overrightarrow{H}(\lambda(t))) = 0.
	\]
	The minimal control for the curve $\gamma = \pi \circ \lambda$ is given by $u_i(t) = \langle \lambda(t), X_k(\pi(\lambda)) \rangle$ and therefore
	\[
	\dfrac{1}{2} \| \dot{\gamma}(t) \|^2 = \dfrac{1}{2} \sum_{i = 1}^m u_i(t)^2 = \dfrac{1}{2} \sum_{i = 1}^m \langle \lambda(t), X_k(\pi(\lambda)) \rangle^2 = H(\lambda(t)).
	\]
\end{proof}

In view of this result, we observe that, contrary to the Riemannian case, there might exist initial covectors $\lambda_0 \in \T_q^*(M)$ such that the corresponding normal geodesic is trivial. This can happen if $\lambda_0 \in H_q^{-1}(0)$.

Since the normal geodesic $\gamma$ has constant-speed, we have
\[
\dot{\gamma}(t) = \dfrac{\diff}{\diff t} \exp_q(t \lambda_0) = \diff_{t \lambda_0} \exp_q(\lambda_0) = \diff_{t \lambda_0} \exp_q(\mathrm{\dot{r}}_{q, \lambda_0}(t))
\]
which is non-zero as long as $\lambda_0 \in \A_q \setminus H_q^{-1}(0)$.

This proves a cotangent version of the first condition of Warner.

\begin{theorem}[Constant speed property]
	The map $\mathrm{exp}_q$ is $\mathcal{C}^\infty$ on $\A_q$ and, for all $\lambda_0 \in \A_q \setminus H_q^{-1}(0)$ and all $t \in I_{q, \lambda_0}$, we have $\mathrm{d}_{t \lambda_0} \mathrm{exp}_q(\mathrm{\dot{r}}_{q, \lambda_0}(t)) \neq 0_{\mathrm{exp}_q(t \lambda_0)}$.
\end{theorem}

The set $\A_q \setminus H^{-1}(0)$ is open and \textit{radially convex} in the following sense: a subset $\mathcal{V} \subseteq \mathrm{T}^*_q(M)$ is radially convex if for every $\lambda_0 \in \mathcal{V}$, we have that $\left\{ t \in \mathds{R} \mid t \lambda_0 \in \mathcal{V} \right\}$ is an open interval.					

\subsection{Jacobi fields and the regularity property} \label{jacobifields}

As will be shown, this property, called \textit{regularity property}, is a feature of \textit{Jacobi fields}, which theory we outline here in the sub-Riemannian context.

Let $\gamma : \interval{0}{T} \to M$ be a normal geodesic and $\lambda(t)$ be its cotangent lift. We can write $\gamma(t) = \exp_q(t \lambda_0)$ for some initial covector $\lambda_0 \in \T_q^*(M)$. Consider a variation of $\gamma(t)$ through normal geodesics
\[
\Gamma(t, s) = \exp_{\sigma(s)}(t V(s))
\]
where $\Lambda(s) = (\sigma(s), V(s))$ is a curve in $\T^*(M)$ with $\Lambda(0) = (q, \lambda_0)$. The curve $\Lambda$ is well defined on a small interval $\ointerval{-\epsilon}{\epsilon}$. A \textit{sub-Riemannian Jacobi field} $J$ along the normal geodesic $\gamma$ can be seen as the variation field of a variation $\gamma$ through normal geodesics:
\[
J(t) = \dfrac{\partial}{\partial s} \exp_{\sigma(s)}(t V(s)) \Big|_{s = 0} = \dfrac{\partial}{\partial s} \exp(\sigma(s), t V(s)) \Big|_{s = 0}.
\]
Remembering that $\exp_q(t v) = \pi \circ \me^{t \overrightarrow{H}} (q, v)$, we have the equalities
\[
J(t) = \dfrac{\partial}{\partial s} \pi\left(\me^{t \overrightarrow{H}} (\sigma(s), V(s)) \right) \Big|_{s = 0} = \dfrac{\partial}{\partial s} \pi\left(\me^{t \overrightarrow{H}} (\Lambda(s)) \right) \Big|_{s = 0} = \diff_{\lambda(t)}\pi\left(\diff_{\lambda_0}\me^{t \overrightarrow{H}} \dot{\Lambda}(0) \right).
\]
The Jacobi field $J$ along $\gamma$ is therefore uniquely determined by its \textit{initial value} $\dot{\Lambda}(0) \in \mathrm{T}_{(q, \lambda_0)}(\mathrm{T}^*(M))$. This implies that the space of \textit{Jacobi fields along the geodesic} $\gamma$, which we denote by $\mathscr{J}(\gamma)$, is a vector space of dimension $2 n$.

On the other hand, the space of \textit{Jacobi fields along the extremal} $\lambda$, denoted this time by $\mathscr{J}(\lambda)$ is the collection of vector fields along $\lambda$ of the form
\[
\mathcal{J}(t) := \diff_{\lambda_0} \me^{t \overrightarrow{H}} \dot{\Lambda}(0),
\]
also uniquely determined by $\dot{\Lambda}(0) \in \mathrm{T}_{(q, \lambda_0)}(\mathrm{T}^*(M))$. The space $\mathscr{J}(\gamma)$ is linearly isomorphic to $\mathscr{J}(\lambda)$ through the pushforward of the bundle projection $\pi : \T^*(M) \to M$. Equivalently, a vector field $\mathcal{J}$ is a Jacobi field along the extremal $\lambda$ if it satisfies
\begin{equation}
	\label{lieeq}
	\dot{\mathcal{J}} := \mathcal{L}_{\overrightarrow{H}} \mathcal{J} = 0,
\end{equation}
where $\mathcal{L}_{\overrightarrow{H}} \mathcal{J}$ is the Lie derivative of a vector field along $\lambda$ in the direction of $\overrightarrow{H}$:
\[
\mathcal{L}_{\overrightarrow{H}} \mathcal{J} (t) = \lim_{\epsilon \to 0} \dfrac{(\diff_{\lambda(t + \epsilon)} \me^{-\epsilon \overrightarrow{H}})[\mathcal{J}(t + \epsilon)] - \mathcal{J}(t)}{\epsilon} = \frac{\diff}{\diff\epsilon}\Big|_{\epsilon=0} (\diff_{\lambda(t + \epsilon)} \me^{-\epsilon \overrightarrow{H}})[\mathcal{J} (t + \epsilon)].
\]

The equation \eqref{lieeq} can be rewritten using the symplectic structure of $\T^*(M)$ and moving frame generalizing, in a \textit{non-canonical} way, the Riemannian parallel transport (see also \cite{subriemjacobi} for more details).

\begin{theorem}
	\label{movingframe}
	Let $\gamma : \interval{0}{T} \to M$ be a normal geodesic and $\lambda(t)$ its cotangent lift. There exists a frame $E_1(t), \dots, E_n(t), F_1(t), \dots, F_n(t)$ along $\lambda(t)$ such that
	\begin{enumerate}[label=\normalfont(\roman*)]
		\item \label{vertical} $\mathrm{Ver}_{\lambda(t)} = \vspan\left\{ E_1(t), \dots, E_n(t) \right\}$, where $\mathrm{Ver}_{\lambda(t)} := \Kern (\diff_{\lambda(t)} \pi) \subseteq \T_{\lambda(t)}(\T^*(M))$ is the \emph{vertical} subspace along $\lambda$;
		\item It is a \emph{symplectic moving frame}:
		\[
		\omega(E_i, E_j) = 0, \qquad \omega(F_i, F_j) = 0, \qquad \omega(E_i, F_j) = \delta_{i,j}.
		\]
	\end{enumerate}
	Furthermore, given such a moving frame, a vector field $\mathcal{J}(t) = \sum_{i = 1}^{n} p_i(t) E_i(t) + x_i(t) F_i(t)$ is a Jacobi field along $\lambda$ if and only the following \emph{Hamilton's equation for Jacobi fields} is satisfied:
	\begin{equation}
		\label{hamiltoneqjacobifields}
		\begin{pmatrix}
			\dot{p} \\
			\dot{x}
		\end{pmatrix}
		=
		\begin{pmatrix}
			- A(t)^T & R(t) \\
			B(t) & A(t)
		\end{pmatrix}
		\begin{pmatrix}
			p \\
			x
		\end{pmatrix}
	\end{equation}
	for some matrix $A(t)$ with $\mathrm{rank} A(t) = \dim \mathcal{D}_{\gamma(t)}$ and some symmetric matrices $B(t), R(t)$.
\end{theorem}

Let $\gamma : \interval{0}{T} \to M$ be a normal geodesic, its cotangent lift $\lambda(t)$, and fix a symplectic moving frame $E_1(t), \dots, E_n(t), F_1(t), \dots, F_n(t)$ along $\lambda$ such as given by \cref{movingframe}. The function $t \mapsto (p(t), x(t))$ from \eqref{hamiltoneqjacobifields} will be called the \textit{coordinates of the Jacobi field} $\mathcal{J}(t)$ along $\lambda$ (resp. $J(t)$ along $\gamma$). The family $\diff_{\lambda(t)}\pi\left(F_1(t)\right), \dots, \diff_{\lambda(t)}\pi\left(F_n(t)\right)$ forms a basis for $\T_{\gamma(t)}(M)$. The \textit{scalar product} $\langle \cdot, \cdot \rangle_{\gamma(t)} : \T_{\gamma(t)}(M) \times \T_{\gamma(t)}(M) \to \R$ designates the positive quadratic form that turns this family into an orthonormal basis along $\gamma$. Furthermore, it coincides with $\langle \cdot, \cdot \rangle_{\D_{\gamma(t)}}$ on $\D_{\gamma(t)}$. All these definitions, including the next one, are dependant on the choice of symplectic moving frame along $\lambda$.

\begin{definition}
	The \textit{derivative} of a Jacobi field $J$ along $\gamma$ with coordinates $(p, x)$ is defined as
	\[
	\nabla J(t) = \sum_{k = 1}^n p_k(t) \diff_{\lambda(t)}\pi\left(F_k(\lambda(t))\right).
	\]  
\end{definition}

The next two results are at the heart of the second condition of Warner.

\begin{lemma}
	\label{orthjacobi}
	Let $\gamma : \interval{0}{T} \to M$ be a normal geodesic. Suppose that $J$ and $\overline{J}$ are two Jacobi fields along $\gamma$. Once we fix a symplectic moving frame $E_1(t), \dots, E_n(t)$, $F_1(t), \dots, F_n(t)$ along the cotangent lift $\lambda(t)$ as given by \cref{movingframe}, we have that
	\[
	\langle \nabla J(t), \overline{J}(t) \rangle_{\gamma(t)} - \langle J(t), \nabla \overline{J}(t) \rangle_{\gamma(t)}
	\]
	is constant in $t \in \interval{0}{T}$.
\end{lemma}

\begin{remark}
	This result is a generalisation of a well-known fact in Riemannian geometry: If $J_1$ and $J_2$ be two \textit{Riemannian} Jacobi fields along a geodesic $\gamma$, then
	$\langle D_t J_1, J_2 \rangle - \langle J_1, D_t J_2 \rangle$ is constant along $\gamma$, where $D_t$ stands for the covariant derivative.
\end{remark}

\begin{proof}
	Let $(p, x)$ (resp. $(\overline{p}, \overline{x})$) be the coordinates of the Jacobi field $J$ (resp. $\overline{J}$) with respect to the moving frame. Hamilton's equation for Jacobi fields \eqref{hamiltoneqjacobifields} states that $\dot{p}(t) = -A(t)^T p(t) + R(t) x(t)$ and $\dot{x}(t) = B(t) p(t) + A(t) x(t)$ and since $R(t)$ and $B(t)$ are symmetric, we have
	\begin{align*}
		\dfrac{\diff}{\diff t} \left[\langle p(t), \overline{x}(t) \rangle_{\R^n} - \langle \overline{p}(t), x(t) \rangle_{\R^n}\right] &= \langle x(t) , R(t) \overline{x}(t) \rangle_{\R^n} - \langle R(t) x(t) , \overline{x}(t) \rangle_{\R^n} \\
		& \quad + \langle p(t) , B(t) \overline{p}(t) \rangle_{\R^n} - \langle B(t) p(t) , \overline{p}(t) \rangle_{\R^n} \\
		& \quad + \langle A(t)^T \overline{p}(t) , x(t) \rangle_{\R^n} - \langle \overline{p}(t) , A(t) x(t) \rangle_{\R^n} \\
		& \quad + \langle p(t) , A(t) \overline{x}(t) \rangle_{\R^n} - \langle A^T p(t) , \overline{x}(t) \rangle_{\R^n} = 0.
	\end{align*}
	This holds for every $t \in \interval{0}{T}$ which concludes the proof.
\end{proof}

Let $\gamma : I \to M$ be a normal geodesic and $a, b \in I$. We use the notation $\mathscr{J}_{a}(\gamma)$ for the vector space of Jacobi fields along $\gamma$ vanishing at time $t = a$ and $\mathscr{J}_{a, b}(\gamma)$ for the subspace of $\mathscr{J}_{a}(\gamma)$ of Jacobi fields along $\gamma : I \to M$ vanishing at both $t = a$ and $t = b$.

\begin{proposition}
	\label{orthodec}
	Let $\gamma : \interval{0}{T} \to M$ be a normal geodesic with initial covector $\lambda_0 \in \T^*_q(M)$ and such that $\gamma(0) = q \in M$. Fix also a symplectic moving frame as provided by \cref{movingframe}. Then, the sets $A_{\gamma(t)} := \left\{ J(t) \ \middle| \ J \in \mathscr{J}_0(\gamma) \right\}$	and $B_{\gamma(t)} := \left\{ \nabla J(t) \ \middle| \ J \in \mathscr{J}_{0, t}(\gamma) \right\}$ are orthogonal complements (with respect to $\langle \cdot, \cdot \rangle_{\gamma(t)}$) in $\T_{\gamma(t)}(M)$, i.e. $\T_{\gamma(t)}(M) = A_{\gamma(t)} \oplus B_{\gamma(t)}$ and $A_{\gamma(t)} = B_{\gamma(t)}^\perp$.
\end{proposition}

\begin{proof}
	Choose a basis $(J_1, \dots, J_k)$ of $\mathscr{J}_{0, t}(\gamma)$ and complete it into a basis for $\mathscr{J}_{0}(\gamma)$ with $(\overline{J}_1, \dots, \overline{J}_{n - k})$.
	
	The family $(\nabla J_1(t), \dots, \nabla J_k(t))$ is linearly independent. Indeed, assume this is not the case, then there exists $a_1, \dots, a_k \in \R$ such that $\sum_{i = 1}^k a_i \nabla J_i(t) = 0$. The Jacobi field $J := \sum_{i = 1}^k a_i J_i$ satisfies $J(0) = 0$, $J(t) = 0$ and $\nabla J(t) = 0$. By Hamilton's equation for Jacobi fields, we can conclude that $J$ is identically zero and therefore $a_1 = \dots = a_k = 0$.
	
	Similarly, the family $(\overline{J}_1, \dots, \overline{J}_{n - k})$ is linearly independent. If this is not the case, then there exists $b_1, \dots, b_{n - k}$ such that $\sum_{i = 1}^{n - k} b_i \overline{J}_i(t) = 0$. Let $\overline{J} = \sum_{i = 1}^{n - k} b_i \overline{J}_i$. We have $\overline{J}(0) = 0$ and $\overline{J}(t) = 0$. So, the Jacobi field $\overline{J} \in \mathscr{J}_{0, t}(\gamma)$ and hence $b_1 = \dots = b_{n - k} = 0$.	
	
	Finally, \cref{orthjacobi} implies that
	\[
	\langle \nabla J_i(t), \overline{J}_j(t) \rangle_{\gamma(t)} - \langle J_i(t), \nabla \overline{J}_j(t) \rangle_{\gamma(t)} = \langle \nabla J_i(0), \overline{J}_j(0) \rangle_{\gamma(0)} - \langle J_i(0), \nabla \overline{J}_j(0) \rangle_{\gamma(0)} = 0
	\]
	and hence $\langle J_i(t), \nabla \overline{J}_j(t) \rangle_{\gamma(t)} = \langle \nabla J_i(t), \overline{J}_j(t) \rangle_{\gamma(t)} = 0$ since $J_1, \dots, J_k \in \mathscr{J}_{0, t}(\gamma)$ by construction. Therefore, we have $A_{\gamma(t)} = B_{\gamma(t)}^\perp$, $A_{\gamma(t)} \cap B_{\gamma(t)} = \{ 0 \}$ and $\T_{\gamma(t)}(M) = A_{\gamma(t)} \oplus B_{\gamma(t)}$.
\end{proof}

It remains to link Jacobi fields and the exponential map. The following two results are analogous to the Riemannian context. Firstly we examine Jacobi fields along $\gamma$ vanishing at its initial time.

\begin{proposition}
	\label{jacobivan}
	Let $\gamma : \interval{0}{T} \to M$ be a normal geodesic with initial covector $\lambda_0 \in \T^*_q(M)$ and such that $\gamma(0) = q \in M$. For every $w \in \T^*_q(M)$, the unique Jacobi field along $\gamma$ with initial value $(0, w) \in \T_q(M) \oplus \T^*_q(M) \cong \mathrm{T}_{(q, \lambda_0)}(\mathrm{T}^*(M))$ is given by
	\[
	J(t) = \diff_{t \lambda_0} \exp_q(t w)
	\]
	where we view $t w$ as an element of $\T_{t \lambda_0}(\T_q^*(M)) \cong \T^*_q(M)$.
\end{proposition}

\begin{proof}
	Consider now a variation of a normal geodesic $\gamma(t) = \mathrm{exp}_q(t \lambda_0)$
	\[
	\Gamma(t, s) = \mathrm{exp}_{\sigma(s)}(t V(s))
	\]
	where $V(s) \in \mathrm{T}^*_{\gamma(s)}(M)$ is a covector field along $\gamma$ satisfying $V(0) = \lambda_0$ and $\dot{V}(0) = w$ and $\sigma(s) = q \in M$ is constant. The variation field of $\Gamma$ is a Jacobi field along $\gamma$
	\[
	\dfrac{\partial}{\partial s} \mathrm{exp}_{q}(t V(s)) \Big|_{s = 0} = \mathrm{d}_{t \lambda_0} \mathrm{exp}_{q}(t w)
	\]
	uniquely determined by its initial value $(\dot{\sigma}(s), \dot{V}(s)) = (0, w)$.
\end{proof}

Secondly, we look at the singularities of the exponential map. These are covectors $\lambda_0 \in \A_q$ such that $\Kern(\diff_{\lambda_0} \exp_q)$ is not trivial. These covectors are called conjugate (co)vectors and the point $\exp_q(\lambda_0)$ is said to be conjugate to $p$.

\begin{proposition}
    \label{singexpp}
	Let $\gamma : \interval{0}{1} \to M$ be a normal geodesic with initial covector $\lambda_0 \in \A_q$ and such that $\gamma(0) = q \in M$. The covector $\lambda_0$ is a critical point for $\exp_q$ if and only if there exists a non trivial Jacobi field $J$ along $\gamma$ such that $J(0) = 0$ and $J(1) = 0$.
\end{proposition}

\begin{proof}
	If $\lambda_0$ is a singularity of $\exp_q$, there exists a vector $\overline{\lambda}_0$ such that $\diff_{\lambda_0} \exp_q(\overline{\lambda}_0) = 0$. In that case, from the previous proposition, the vector field
	\[
	J(t) = \diff_{t \lambda_0} \exp_q(t \overline{\lambda}_0)
	\]
	is a non trivial Jacobi field field such that $J(0) = 0$ and $J(1) = 0$. The converse implication is similar.
\end{proof}

For $A \in \Kern(\diff_{\lambda_0} \exp_q)$, we denote by $J_A(t) \in \mathscr{J}_{0, 1}(\gamma)$ the Jacobi field along $\gamma$ with initial value $(0, A)$. We finish this section by proving the following cotangent version of Warner's second condition.

\begin{proposition}[Regularity property]
	\label{propregularity}
	For every $\lambda_0 \in \A_q \setminus \{ 0 \}$ and every symplectic moving frame along the cotangent lift $\lambda(t)$ of the normal geodesic $\gamma(t):= \exp_q(t \lambda_0)$, as given by \cref{movingframe}, the map 
	\[
	\Kern(\diff_{\lambda_{0}} \exp_q) \to \T_{\exp_q(\lambda_{0})}(M)/\diff_{\lambda_{0}} \exp_q(\T_{\lambda_{0}}(\T^*_q(M)))
	\]
	sending $A$ to $\nabla J_A(1) + \diff_{\lambda_{0}} \exp_q(\T_{\lambda_{0}}(\T^*_q(M)))$, is a linear isomorphism.
\end{proposition}

\begin{proof}
	Let $\lambda_0 \in \T_q^*(M)$ and $(A_1, \dots, A_k)$ be a basis for $\Kern(\diff_{\lambda_{0}} \exp_q)$. We can view them as elements of $\T_q^*(M)$ via the identification $\T_{\lambda_{0}}(\T_q^*(M)) \cong \T_q^*(M)$ when necessary. 
	
	For $i = 1, \dots, k$, the Jacobi fields
	\[
	J_i(t) := \diff_{t \lambda_0} \exp_q(t A_i)
	\]
	vanish at their initial time and have the the initial value $(0, A_i) \in \T_q(M) \oplus \T^*_q(M) \cong \mathrm{T}_{(q, \lambda_0)}(\mathrm{T}^*(M))$ (see \cref{jacobivan}). They also vanish at the final time $t = 1$ since $A_1, \dots, A_k \in \Kern(\diff_{\lambda_{0}} \exp_q)$. Therefore, $J_i \in \mathscr{J}_{0, 1}(\gamma)$ for every $i = 1, \dots, k$. Using a moving frame along $\lambda$ given by \cref{movingframe}, we define a linear map
	\[
	\theta : \Kern(\diff_{\lambda_{0}} \exp_q) \to \T_{\exp_q(v)}(M)/\diff_{\lambda_{0}} \exp_q(\T_{\lambda_{0}}(\T^*_q(M)))
	\]
	via $\theta(A_i) := \nabla J_i(1) + \diff_{\lambda_{0}} \exp_q(\T_{\lambda_{0}}(\T^*_q(M)))$. \cref{orthodec} implies that the family $(\nabla J_1(1), \dots, \nabla J_k(1))$ is linearly independent and that $\Kern(\diff_{\lambda_{0}} \exp_q) \cong B_{\gamma(1)}$ by the orthogonal decomposition. Therefore, the map $\theta$ is a linear isomorphism.
\end{proof}

\begin{remark}
	In light of \cref{orthodec}, \cref{singexpp}, \cref{propregularity} and its proof, we can see that
	\[
	A_{\gamma(1)} = \diff_{\lambda_{0}} \exp_q(\T_{\lambda_{0}}(\T^*_q(M)))
	\]
	and 
	\begin{align*}
		B_{\gamma(1)} & \cong \left\{ \sum_{k = 1}^n p_k(1) E_k(1) \mid (p(t), x(t)) \text{ satisfies \eqref{hamiltoneqjacobifields} with } x(0) = x(1) = 0 \right\} \\
		& = \Kern(\diff_{\lambda_0} \exp_q).
	\end{align*}
	In particular, the subspace $B_{\gamma(1)}$ does not depend on the choice of moving frame along $\lambda(t)$.
\end{remark}

Let us describe the results of this section in a different way. The introduction of a moving frame as given in \cref{movingframe} is equivalent to choosing a horizontal complement of some Riemannian metric $g$ on $M$. Indeed, given a moving frame $E_1(t), \dots, E_n(t), F_1(t), \dots, F_n(t)$, we can naturally obtain a Riemannian metric $g$ that extends $\langle \cdot, \cdot \rangle_{\gamma(t)}$ along the curve $\gamma(t)$ to the whole manifold $M$. Conversely, choosing a Riemannian metric $g$ on $M$ induces an isomorphism:
\[
\sharp : \mathrm{Ver}_{\lambda} \to \mathrm{T}_{q}(M) : \xi \mapsto \xi^\sharp, \qquad p := \pi(\lambda), \ \lambda \in \mathrm{T}^*(M)
\]
where $\xi^\sharp$ is the unique element of $\mathrm{T}_q(M)$ such that $g(\xi^\sharp, X) = \xi(X)$, for every $X \in \mathrm{T}_q(M)$, the spaces $\mathrm{Ver}_\lambda$ and $\mathrm{T}^*_q(M)$ being canonically identified. Now, if $X \in \mathrm{T}_q(M)$ and if $F \in \mathrm{T}_\lambda(\mathrm{T}^*(M))$ such that $\mathrm{d}_\lambda \pi (F) = X$, one can show that
\begin{equation}
	\label{idgomega}
	g(\xi^\sharp, X) = \xi(X) = \omega(\xi, F).
\end{equation}
Therefore, from an orthonormal basis $X_1, \dots, X_n$ along $\gamma(t)$, we can construct the moving frame $E_1(t), \dots, E_n(t), F_1(t), \dots, F_n(t)$ by setting $E_i^\sharp(t) = X_i(t)$ and $X_i(t) = \mathrm{d}_{\lambda(t)}\pi(F_i(t))$, and \cref{movingframe} would follow from \cref{idgomega}. Different choices of $g$ therefore corresponds to different choices of Darboux moving frames.

When such a Riemannian metric $g$ on $M$ has been fixed, the derivative of a Jacobi field along $\gamma$ with initial value $\xi_0$ corresponds to
\[
\nabla J(t) = \left(\mathrm{d}_{\lambda_0} \mathrm{e}^{t \overrightarrow{H}} [\xi_0]\right)^\sharp.
\]
Finally, \cref{orthodec} is saying that the map
\[
\mathrm{Ker}(\mathrm{d}_{\lambda_0} \mathrm{exp}_q(\lambda_0)) \to \mathrm{T}_{\mathrm{exp}_q(\lambda_0)}(M) : \xi_0 \mapsto \left(\mathrm{d}_{\lambda_0} \mathrm{e}^{t \overrightarrow{H}} [\xi_0]\right)^\sharp
\]
has its image $g$-perpendicular to $\mathrm{Im}(\mathrm{d}_{\lambda_0} \mathrm{exp}_q)$, leading to \cref{propregularity}				
																
\subsection{Maslov index and the continuity property} \label{R3}

We now approach Jacobi fields and conjugate points via Lagrange Grassmannian geometry. For the definitions and properties related to Jacobi curves, we refer the reader to \cite[Chapter 15]{comprehensive2020}, while the Maslov index, in its full generality, is developed in \cite[Chapter 5]{piccione}, for example.

We start with Lagrange Grassmannian geometry.

Let $(\Sigma, \omega)$ be a symplectic vector space of dimension $2n$. A Lagrangian subspace of $\Sigma$ is a subspace of dimension $n$ on which the restriction of $\omega$ vanishes. The collection of all Lagrangian subspaces of $\Sigma$ is called the Lagrange Grassmannian of $\Sigma$, and is denoted by $L(\Sigma)$. It is a compact manifold of dimension $n(n+1)/2$. 
Furthermore, if $\Lambda \in L(\Sigma)$, there is a linear isomorphism between the tangent space $\T_\Lambda(L(\Sigma))$ and the space of bilinear forms $Q(\Lambda)$ on $\Lambda$. The linear isomorphism is given by 
\[
\T_\Lambda(L(\Sigma)) \to Q(\Lambda) : \dot{\Lambda}(0) \mapsto \underline{\dot{\Lambda}},
\]
where $\underline{\dot{\Lambda}}(z) := \omega(z(0), \dot{z}(0))$, where we consider a smooth curve $\Lambda(t)$ in $L(\Sigma)$ such that $\Lambda(0) = \Lambda$ and a smooth extension $z(t) \in \Lambda(t)$ such that $z(0) = 0$. It can be shown that $\underline{\dot{\Lambda}}$ is a well-defined quadratic map, independent of the extension considered.

For $k = 0, \dots, n$, we define the following subsets of $L(\Sigma)$
\[
\Lambda^k(L_0) = \{ L \in L(\Sigma) \mid \dim(L \cap L_0) = k\} \text{ and } \Lambda^{\geqslant k}(L_0) = \bigcup_{i = k}^n \Lambda^i(L_0).
\]
For each $k = 0, \cdots, n$, the spaces $\Lambda^k(L_0)$ are embedded manifolds of $L(\Sigma)$ with codimension $\frac{1}{2}k(k+1)$. Their tangent space is given by
\[
T_L(\Lambda^k(L_0)) = \{B \in B_{\mathrm{sym}}(L) \mid B_{(L_0 \cap L) \times (L_0 \cap L)} = 0\},
\]
for all $L \in \Lambda^k(L_0)$, where $B_{\mathrm{sym}}(L)$ is the space of symmetric forms on $L$.

Let $\Lambda(\cdot) : \interval{a}{b} \to L(\Sigma)$ be a curve of class $\mathcal{C}^1$. We say that \textit{$\Lambda(\cdot)$ intercepts $\Lambda^{\geqslant 1}(L_0)$ transversally at the instant $t = t_0$} if $\Lambda(t_0) \in \Lambda^1(L_0)$ and if the symmetric bilinear form $\underline{\dot{\Lambda}}(t_0)$ is non zero in the space $\Lambda(t_0) \cap L_0$. This intersection is \textit{positive} (resp. \textit{negative}, \textit{non-degenerate}) if $\underline{\dot{\Lambda}}(t_0)$ is positive definite (resp. negative definite, \textit{non-degenerate}).

We recall that the first group of relative homology of the pair $(L(\Sigma), \Lambda^0(L_0))$ is an infinite cyclic group (\cite[Corrolary 1.5.3.]{piccione}) and is denoted by $H_1(L(\Sigma), \Lambda^0(L_0))$.

\begin{definition}
	Let $L_0 \in L(\Sigma)$. We say that a curve $\Lambda(\cdot) : \interval{a}{b} \to L(\Sigma)$ of class $\mathcal{C}^1$ with endpoints in $\Lambda^0(L_0)$ is a \textit{positive generator} of $H_1(L(\Sigma), \Lambda^0(L_0))$ if $\Lambda(\cdot)$ transversally and positively intercepts $\Lambda^{\geqslant 1}(L_0)$ only once.
\end{definition}

Positive generators are homologous in $H_1(L(\Sigma), \Lambda^0(L_0))$ and thus any of these curves defines a generator of $H_1(L(\Sigma), \Lambda^0(L_0)) \cong \mathds{Z}$ (\cite[Lemma 5.1.11]{piccione}). This fact ensures us that the next object is well defined.

\begin{definition}
	An isomorphism
	\begin{equation}
		\label{maslov}
		\mu_{L_0} : H_1(L(\Sigma), \Lambda^0(L_0)) \to \Z
	\end{equation}
	is defined by requiring that any positive generator of $H_1(L(\Sigma), \Lambda^0(L_0))$ is sent to $1 \in \Z$. If $\Lambda(\cdot) : \interval{a}{b} \to L(\Sigma)$ is a continuous curve with endpoints in $\Lambda^0(L_0)$, we denote by $\mu_{L_0}(\Lambda(\cdot)) \in \Z$ the integer number that corresponds to the homology class of $\Lambda(\cdot)$ by the isomorphism \eqref{maslov}. The number $\mu_{L_0}(\Lambda(\cdot))$ is called the \textit{Maslov index} of the curve $\Lambda(\cdot)$ relative to the Lagrangian $L_0$.
\end{definition}

The key properties of the Maslov index, including the fact that it is homotopy invariant, is summarized in the next theorem.

\begin{theorem}[{{\cite[Lemma 5.1.13 and Corollary 5.1.18]{piccione}}}]
	\label{propmaslov}
	Let $\Lambda : \interval{a}{b} \to L(\Sigma)$ be a curve with endpoints in $\Lambda^0(L_0)$. We have
	\begin{theoenum}
		\item If $\sigma : \interval{c}{d} \to \interval{a}{b}$ is a continuous map with $\sigma(c) = a$ and $\sigma(d) = b$, then $\mu_{L_0}(\Lambda(\sigma(\cdot))) = \mu_{L_0}(\Lambda(\cdot))$;
		\item If $\Lambda'(\cdot) : \interval{c}{d} \to L(\Sigma)$ is a curve with endpoints in $\Lambda^0(L_0)$ and such that $\Lambda(b) = \Lambda'(c)$, then $\mu_{L_0}((\Lambda \star \Lambda') (\cdot)) = \mu_{L_0}(\Lambda(\cdot)) + \mu_{L_0}(\Lambda'(\cdot))$, where $\star$ denotes the concatenation of curves;
		\item $\mu_{L_0}(\Lambda(\cdot)^{-1}) = - \mu_{L_0}(\Lambda(\cdot))$ where $\cdot^{-1}$ denotes the reversed curve;
		\item If $\Ima(\Lambda(\cdot)) \subset \Lambda^0(L_0)$, then $\mu_{L_0}(\Lambda(\cdot)) = 0$;
		\item \label{5} If $\Lambda'(\cdot) : \interval{a}{b} \to L(\Sigma)$ is a curve homotopic to $\Lambda(\cdot)$ with free endpoints in $\Lambda^0(L_0)$, i.e. there is exists a continuous function $H : \interval{0}{1} \times \interval{a}{b} \to L(\Sigma)$ such that $H(0, t) = \Lambda(t)$, $H(1, t) = \Lambda'(t)$ for every $t \in \interval{a}{b}$ and $H(s, a), H(s, b) \in \Lambda^0(L_0)$ for every $s \in \interval{0}{1}$, then $\mu_{L_0}(\Lambda'(\cdot)) = \mu_{L_0}(\Lambda(\cdot))$;
		\item We have $\mu_{L_0}(\Lambda(\cdot)) = \mu_{S(L_0)}(S \circ \Lambda(\cdot))$ if $S : (\Sigma, \omega) \to (\Sigma', \omega')$ is a symplectomorphism;
		\item if $\Lambda(\cdot) : \interval{a}{b} \to L(\Sigma)$ is of class $\mathcal{C}^1$ with endpoints in $\Lambda^0(L_0)$ that has only non-degenerate intersection with $\Lambda^{\geqslant 1}(L_0)$ and if $\Lambda(t) \in \Lambda^{\geqslant 1}(L_0)$ only at a finite number of $t \in \interval{a}{b}$, then
		\[
		\mu_{L_0}(\Lambda(\cdot)) = \sum_{t \in \interval{a}{b}} \mathrm{sgn}(\underline{\dot{\Lambda}}(t)|_{(L_0 \cap \Lambda(t)) \times (L_0 \cap \Lambda(t))}),
		\]
		where $\mathrm{sgn}(B)$ is the \textit{signature} of $B$, that is to say, $\mathrm{sgn}(B) = \eta_+(B) - \eta_-(B)$ with
		\[
		\eta_{+/-}(B) = \sup\{ \dim(W) | B_{W \times W} \text{ is negative/positive definite}\}.
		\]
	\end{theoenum}
\end{theorem}

\begin{remark}
    \cref{5} states the homotopy invariance of the Maslov index. This property was firstly proved in \cite{arnold1997}.
\end{remark}

We now turn our attention to the concept of Jacobi curve and the continuity property. As seen in the previous section, a conjugate vector occurs when there is a non trivial Jacobi field along the geodesic that vanishes both at its initial and endpoint. It seems therefore natural to study the evolution of the space of all Jacobi fields at a time $t$ that vanish at its initial time:
\begin{equation}
\label{curve}
\mathrm{L}_{(q, \lambda_0)}(t) := \{ \mathcal{J}(t) \mid \mathcal{J} \text{ is a Jacobi field along } \lambda(t) := \me^{t \overrightarrow{H}}(q, \lambda_0) \text{ and } J(0) = 0 \},
\end{equation}
for $q \in M$ and $\lambda_0 \in \mathscr{A}_q$. For every $t \in \interval{0}{1}$, the set $\mathrm{L}_{(q, \lambda_0)}(t)$ is a Lagrangian subspace of $\T_{\lambda(t)}(\T^*(M))$. In order to be able to use the geometry and theory of Lagrangian Grassmannian, we will work with an alternative curve that lives in the fixed Grassmannian $\T_{(q, \lambda_0)}(\T^*(M))$.

\begin{definition}
	Let $\gamma : \interval{0}{T} \to M$ be a normal geodesic with initial covector $\lambda_0 \in \T^*_q(M)$. The Jacobi \textit{curve} along the cotangent lift $\lambda : \interval{0}{T} \to \T^*(M)$ is defined by
	\[
	\mathrm{J}_{(q, \lambda_0)}(t) := \diff_{\lambda(t)} \me^{-t \overrightarrow{H}} [\mathrm{Ver}_{\lambda(t)}]
	\]
	where $\mathrm{Ver}_{\lambda} := \T_\lambda(\T^*_{\pi(\lambda)}(M))$ is the vertical subspace of $\T_\lambda(\T^*(M))$.
\end{definition}

\begin{remark}
	In particular, we see that any point of the Jacobi curve $\mathrm{J}_{(q, \lambda_0)}(t)$ is a subspace of $\T_{(q, \lambda_0)}(\T^*(M))$ and that $\mathrm{J}_{(q, \lambda_0)}(0) = \T_{(q, \lambda_0)}(\T^*(M))$. The subscript $(q, \lambda_0)$ on $\mathrm{J}_{(q, \lambda_0)}(t)$ is a reminder of the fact that in order to defined the Jacobi curve, one only needs to specify a point $q \in M$ covector $\lambda_0 \in \T^*_q(M)$, the cotangent lift being $\lambda(t) = \me^{t \overrightarrow{H}}(q, \lambda_0)$ and the normal geodesic $\gamma(t) = \pi(\lambda(t))$. The Jacobi curve can be seen as the evolution of the space of all Jacobi fields at time 0 that vanish at time $t$:
	\[
	\mathrm{J}_{(q, \lambda_0)}(t) = \{ \mathcal{J}(0) \mid \mathcal{J} \text{ is a Jacobi field along } \lambda(s) := \me^{s \overrightarrow{H}}(q, \lambda_0) \text{ and } J(t) = 0 \}.
	\]
\end{remark}

We state basic properties of Jacobi curves.

\begin{proposition}[{\cite[Proposition 15.2.]{comprehensive2020}}]
	The Jacobi curve $\mathrm{J}_{(q, \lambda_0)}(t)$ satisfies the following properties:
	\begin{enumerate}[label=\normalfont(\roman*)]
		\item $\mathrm{J}_{(q, \lambda_0)}(t + s) = \diff_{\lambda(t + s)} \me^{-t \overrightarrow{H}} [\mathrm{J}_{(q, \lambda_0)}(s)]$;
		\item $\underline{\dot{\mathrm{J}}}_{(q, \lambda_0)}(0) = - 2 H_{q}$ as quadratic forms on $\mathrm{Ver}_{(q, \lambda_0)} \cong \T^*_q(M)$;
		\item $\mathrm{rank} \ \underline{\dot{\mathrm{J}}}_{(q, \lambda_0)}(t) = \mathrm{rank} \ H_{q}$,
	\end{enumerate}
	for every $t, s \in \R$ such that both sides of the statements	are defined.
\end{proposition}

The previous result tells us that the Jacobi curve $\mathrm{J}_{(q, \lambda_0)}(t)$ is a \textit{monotone nonincreasing curve}, i.e. $\underline{\dot{\mathrm{J}}}_{(q, \lambda_0)}(t)$ nonpositive quadratic form for every $t \in \interval{0}{T}$.

The relation between the Jacobi curve and the conjugates vectors is given in next proposition.

\begin{proposition}[{\cite[Proposition 15.6.]{comprehensive2020}}]
	\label{jacurveconnj}
	For every $q \in M$, a cotangent vector $s \lambda_0 \in \A_q$ is a conjugate vector if and only if
	\[
	\mathrm{J}_{(q, \lambda_0)}(s) \cap \mathrm{J}_{(q, \lambda_0)}(0) \neq \left\{ 0 \right\}.
	\]
	Furthermore, the order of $s \lambda_0$ as a singularity of $\exp_q$ is equal to $\dim(\mathrm{J}_{(q, \lambda_0)}(s) \cap \mathrm{J}_{(q, \lambda_0)}(0))$.
\end{proposition}

In other words, $s \lambda_0$ is conjugate if and only if the Jacobi curve $\mathrm{J}_{(q, \lambda_0)}(t)$ is in $\Lambda^{\geqslant 1}(L_0)$ for $t = s$. \cref{jacurveconnj} alongside the condition of abnormality for geodesics shows that if we have a segment of points that are conjugate to the initial one, then the segment is also abnormal.

\begin{corollary}[{\cite[Proposition 15.7.]{comprehensive2020}}]
	\label{discrete}
	Let $\mathrm{J}_{(q, \lambda_0)}(t)$ be a Jacobi curve associated with $(q, \lambda_0) \in \T^*_q(M)$ and let $\gamma(t)$ the corresponding normal geodesic. Then, $\gamma|_{\interval{0}{s}}$ is abnormal if and only if $\mathrm{J}_{(q, \lambda_0)}(t) \cap \mathrm{J}_{(q, \lambda_0)}(0) \neq \left\{ 0 \right\}$ for every $t \in \interval{0}{s}$.
\end{corollary}

In particular, a geodesic that does contain an abnormal segment has an infinite number of conjugate points while a strongly normal geodesic has a discrete set of conjugate points. 

The manifold $\Sigma := \T_{(q, \lambda_0)}(\T^*(M))$ has dimension $2 n$ and the cotangent symplectic form on $\T^*(M)$ induces a symplectic bilinear form on $\Sigma$. The Jacobi curve $\mathrm{J}_{(q, \lambda_0)}(t)$ defines a one parameter family of $n$-dimensional subspaces of $\Sigma$. This is because the Hamiltonian flow $\me^{-t \overrightarrow{H}}$ is a symplectic transformation and the vertical space $\mathrm{Ver}_{\lambda}$ is a Lagrangian subspace of $\Sigma$. By consequence, the Jacobi curve $\mathrm{J}_{(q, \lambda_0)}(t)$ is a smooth curve in the Lagrangian Grassmannian $L(\Sigma)$.
 
The Maslov index can therefore be computed for Jacobi curves by applying \cref{propmaslov}, proving the following proposition.

\begin{proposition}
	\label{counting}
	For $\Sigma := \T_{(q, \lambda_0)}(\T^*(M))$ with $\lambda_0 \in \mathscr{A}_q \setminus H _{q}^{-1}(0)$, $L_0 := \mathrm{J}_{(q, \lambda_0)}(r)$ and $\Lambda_{r, s}(t) := \mathrm{J}_{(q, \lambda_0)}|_{\interval{r}{s}}(t)$ where $r < s$ are chosen such that $r \lambda_0$ and $s \lambda_0$ are not conjugate vectors of $\exp_q$, we have
	\[
	\mu_{L_0}(\Lambda_{r, s}(\cdot)) = - \sum_{\substack{t \in \interval{r}{s}\\ \mathrm{Ker}(\diff_{t \lambda_0} \exp_{q}) \neq \left\{ 0 \right\} }} \dim(\Lambda_{r, s}(t) \cap \Lambda_{r, s}(r)).
	\]
\end{proposition}


We are now ready to prove Warner's third condition of regularity.

\begin{proposition}[Continuity property]
	\label{continuity}
	Let $M$ be sub-Riemannian manifold and $\lambda_0 \in \mathscr{A}_q \setminus H_q^{-1}(0)$ a covector such that the corresponding geodesic $\gamma(t) := \exp_q(t \lambda_0)$ is strongly normal. Then, there exists a radially convex neighbourhood $\mathcal{V}$ of $\lambda_0$ such that for every ray $r_{q, \overline{\lambda}_0}$ which intersects $\mathcal{V}$ that does not contain abnormal subsegments in $\mathcal{V}$, the number of  singularities of $\exp_q$ (counted with multiplicities) on $\mathrm{Im}(r_{q, \overline{\lambda}_0}) \cap \mathcal{V}$ is constant and equals the order of $\lambda_0$ as a singularity of $\exp_q$, i.e. $\mathrm{dim} (\mathrm{Ker}(\diff_{\lambda_0} \exp_{q}))$.
\end{proposition}

\begin{remark}
    \cref{continuity} establishes the homotopy of Jacobi curves corresponding to a small path of initial conditions in the cotangent space. The homotopy invariance of the Maslov index was also used in a similar way in \cite[Proposition 2.6.]{sachkov2008}.
\end{remark}

\begin{proof}
    By \cref{discrete}, the conjugate points to $\gamma(0)$ along $\gamma(t) = \exp_q(t \lambda_0)$ are isolated. Thus setting $t_{\pm} := 1 \pm \delta$ for $\delta > 0$ sufficiently small, the Jacobi curve $\mathrm{J}_{(q, \lambda)}|_{\interval{t_-}{t_+}}$  is transverse to $L_0 := \mathrm{J}_{(q, \lambda_0)}(0) = \mathrm{T}_{\lambda_0}(\mathrm{T}^*_q(M))$. By \cref{movingframe}, there exists a symplectic moving frame along $\lambda(t)$ such that a vector field $\mathcal{J}(t) = \sum_{i = 1}^{n} p_i(t) E_i(t) + x_i(t) F_i(t)$ is a Jacobi field if
	\begin{equation}
		\label{jacobmatrix}
	\begin{pmatrix}
		\dot{p} \\
		\dot{x}
	\end{pmatrix}
	=
	\begin{pmatrix}
		- A_{\lambda_0}(t)^T & R_{\lambda_0}(t) \\
		B_{\lambda_0}(t) & A_{\lambda_0}(t)
	\end{pmatrix}
	\begin{pmatrix}
		p \\
		x
	\end{pmatrix}
	\end{equation}
	for some matrix $A_{\lambda_0}(t)$ and some symmetric matrices $B_{\lambda_0}(t)$ and $R_{\lambda_0}(t)$. Such a choice of moving frame can enable us to identify the Jacobi curve $\mathrm{J}_{(q, \lambda_0)}(t)$ with
	\[
	\widetilde{\mathrm{J}}_{(q, \lambda_0)}(t) := \{ (p(0), x(0)) \mid (p, x) \text{ is a solution of \eqref{jacobmatrix} such that } x(t) = 0 \}.
	\]
	With the same moving frame we can identify $\mathrm{L}_{(q, \lambda_0)}(t)$ with
	\[
	\widetilde{\mathrm{L}}_{(q, \lambda_0)}(t) := \{ (p(t), x(t)) \mid (p, x) \text{ is a solution of \eqref{jacobmatrix} such that } x(0) = 0 \}.
	\]
	In these coordinates, the space $\mathrm{J}_{(q, \lambda_0)}(0)$ is therefore identified with the vertical subspace $\mathds{R}^n \times \left\{ 0 \right\}$ of $\mathds{R}^{2n}$. Both $\widetilde{\mathrm{J}}_{(q, \lambda_0)}(\cdot)$ and $\widetilde{\mathrm{L}}_{(q, \lambda_0)}(\cdot)$ are smooth curves in $L(\mathds{R}^{2 n})$ and
	\begin{equation}
	\label{identification}
	\mu_{L_0}(\mathrm{J}_{(q, \lambda_0)}(\cdot)|_{\interval{t_-}{t_+}}) = \mu_{\widetilde{L}_0}(\widetilde{\mathrm{J}}_{(q, \lambda_0)}(\cdot)|_{\interval{t_-}{t_+}}) = \mu_{\widetilde{L}_0}(\widetilde{\mathrm{L}}_{(q, \lambda_0)}(\cdot)|_{\interval{t_-}{t_+}})
	\end{equation}
	with $\widetilde{L}_0 := \mathds{R}^n \times \{0\}$. Indeed, the moving frame $(E_1, \dots, E_n, F_1, \dots, F_n)$ along the lift $\lambda(t)$ induces a symplectomorphism between $\T_{\lambda(t)}(\T^*(M))$ and $\mathds{R}^{2 n}$. \cref{propmaslov} then justifies the first equality in \eqref{identification}. The second equality follows when observing that $\widetilde{\mathrm{J}}_{(q, \lambda_0)}(t) \cap \widetilde{L}_0 \neq \left\{ 0 \right\}$ if and only if $\widetilde{\mathrm{L}}_{(q, \lambda_0)}(t) \cap \widetilde{L}_0 \neq \left\{ 0 \right\}$ and for these conjugate times
	\[
	\dim(\widetilde{\mathrm{J}}_{(q, \lambda_0)}(t) \cap \widetilde{L}_0) = \dim(\widetilde{\mathrm{L}}_{(q, \lambda_0)}(t) \cap \widetilde{L}_0).
	\] 
	
	We then extend the moving frame along the extremal $\lambda(t)$ to a symplectic frame on an open set $\mathcal{U} \subseteq \T^*(M)$ containing $\Ima(\lambda)$. Let us explain why this is possible. If the extremal $\lambda(t)$ does not self intersect, then there is of course such an extension. If an intersection occurs, because $\lambda(t)$ is the solution of a first order ODE on $\T^*(M)$, the extremal $\lambda(t)$ must be a closed curve. In that case, a symplectic extension of the moving frame exists since $\T^*(M)$ naturally oriented.
	
	Now, let $\mathcal{W}$ be a convex neighbourhood of $\lambda_0$ in $\mathrm{T}^*_q(M)$ and for $\overline{\lambda}_0 \in \mathcal{W}$ we set $\overline{\lambda}(t) := \me^{t\overrightarrow{H}}(\overline{\lambda}_0)$, $\overline{\gamma}(t) := \pi(\gamma(t))$, and $l(s) = (1 - s) \lambda_0 + s \overline{\lambda}_0$. The smooth symplectic frame $(E_1(q'), \dots, E_n(q'), F_1(q'), \dots, F_n(q'))$ for $q' \in \mathcal{U}$ is a moving frame along $\lambda(t)$ when $q' = \lambda(t)$ and along $\overline{\lambda}(t)$ when $q' = \overline{\lambda}(t)$, in the sense of \cref{movingframe}. With respect to them, we construct $\mathrm{J}_{(q, \overline{\lambda}_0)}(\cdot)$ and $\overline{L}_0 := \mathrm{J}_{(q, \overline{\lambda}_0)}(0) = \mathrm{T}_{\overline{\lambda}_0}(\mathrm{T}^*_q(M))$ as before and we claim that the curves $\widetilde{\mathrm{L}}_{(q, \lambda_0)}(\cdot)|_{\interval{t_-}{t_+}}$ and $\widetilde{\mathrm{L}}_{(q, \overline{\lambda}_0)}(\cdot)|_{\interval{t_-}{t_+}}$ are homotopic with free endpoints in $\Lambda^0(\widetilde{L}_0)$. Set
	\[
	H : \interval{0}{1} \times \interval{t_-}{t_+} : (s, t) \mapsto \widetilde{\mathrm{L}}_{(q,l(s))}(t)|_{\interval{t_-}{t_+}}.
	\]
	By continuous dependence on $(q, \lambda_0) \in \mathrm{T}^*(M)$ the Jacobi curve $\mathrm{J}_{(q, \lambda_0)}(\cdot)$ has endpoints transverse to $L_0$ for any sufficiently small neighbourhood of $\lambda_0$ in $\mathrm{T}^*_q(M)$. Thus, the curves $\widetilde{\mathrm{L}}_{(q, \lambda_0)}(\cdot)|_{\interval{t_-}{t_+}}$ and $\widetilde{\mathrm{L}}_{(q, \overline{\lambda}_0)}(\cdot)|_{\interval{t_-}{t_+}}$ have endpoints in $\Lambda^0(\widetilde{L}_0)$ since $\gamma(t_-)$, $\overline{\gamma}(t_-)$, $\gamma(t_+)$, $\overline{\gamma}(t_+)$ are not conjugate to $p$, as long as the neighbourhood $\mathcal{W}$ is chosen sufficiently small. Note that we are not saying that $\overline{\gamma}|_{\ointerval{t_-}{t_+}}$ has no abnormal subsegments. Since $l(s) \in \mathcal{W}$ for all $s \in \interval{0}{1}$ by convexity, the same goes for the geodesic starting at $p$ with initial covector $l(s)$. Therefore, the curve $\widetilde{\mathrm{L}}_{(q,l(s))}(\cdot)|_{\interval{t_-}{t_+}}$ has its endpoints in $\Lambda^0(\widetilde{L}_0)$ too. It remains to prove the continuity of $H$. Because the solutions of the geodesic equation depend continuously on the initial covector $\lambda_0$, the matrices in \eqref{jacobmatrix} and therefore $\widetilde{\mathrm{L}}_{(q, \lambda_0)}(\cdot)|_{\interval{t_-}{t_+}}$ also depend continuously on $\lambda_0$. Consequently, $H$ is a homotopy and we can use \cref{propmaslov} and \eqref{identification} to conclude that
	\[
	\mu_{L_0}(\mathrm{J}_{(q, \lambda_0)}(\cdot)|_{\interval{t_-}{t_+}}) = \mu_{\widetilde{L}_0}(\widetilde{\mathrm{L}}_{(q, \lambda_0)}(\cdot)|_{\interval{t_-}{t_+}}) = \mu_{\widetilde{L}_0}(\widetilde{\mathrm{L}}_{(q, \overline{\lambda}_0)}(\cdot)|_{\interval{t_-}{t_+}}) = \mu_{\overline{L}_0}(\mathrm{J}_{(q, \overline{\lambda}_0)}(\cdot)|_{\interval{t_-}{t_+}}).
	\]
	If the ray $r_{q, \overline{\lambda}_0}$ does not contain abnormal subsegments in $\mathcal{W}$, then we have that $\gamma$ and $\overline{\gamma}$ have the same number of conjugate times in $\interval{t_-}{t_+}$, since the Maslov index of the Jacobi curve $\mathrm{J}_{\lambda}(t)$ with endpoints transverse to $\mathrm{J}_{\lambda}(0)$ is equal to the number of intersection (counted with multiplicity, by \cref{counting}), and since such intersections corresponds to conjugate times (the dimension of the intersection being the multiplicity, by \cref{jacurveconnj} and \cref{discrete}). The truncated cone over $\mathcal{W}$ given by $\mathcal{V} := \left\{ \tau \overline{\lambda}_0 \mid \overline{\lambda}_0 \in \mathcal{W}, \tau \in \interval{t_-}{t_+} \right\}$ is the required radially convex set given by the continuity property.
\end{proof}

\begin{remark}
    In fact, even if $\overline{\lambda}_0$ corresponds to a normal geodesic with abnormal subsegments in $\mathcal{V}$, the proof of \cref{continuity} says that $\mathrm{J}_{(q, \overline{\lambda}_0)}|_{\interval{t_-}{t_+}}$ and $\mathrm{J}_{(q, \lambda_0)}|_{\interval{t_-}{t_+}}$ have the same Maslov index with respect to $\mathrm{J}_{(q, \overline{\lambda}_0)}(0)$ and $\mathrm{J}_{(q, \lambda_0)}(0)$.
\end{remark}

This also completes the proof of \cref{regsubriemexp}. It is reasonable to ask whether the cotangent version of Warner's regularity conditions implies a sub-Riemannian analogue of Morse-Littauer's theorem, i.e. the non-injectivity of the sub-Riemannian exponential map on any neighbourhood of a conjugate covector. Warner's approach involves giving the normal forms of the exponential map on neighbourhood of (regular) conjugate vectors. It is not obvious to us that \cref{regsubriemexp} would easily provide such a local description of the sub-Riemannian exponential map about its singularities. However, we are able to pursue Warner's program for a specific example: the three dimensional Heisenberg group.				
									
\section{Exponential map of the Heisenberg group} \label{examples}

In this section, we study the three dimensional Heisenberg group and prove \cref{thm:noninjH}. 

\subsection{Geometry of the Heisenberg group}

The Heisenberg group $\He$ is the sub-Riemannian structure, defined on $\R^3$, that is generated by the global vector fields
\[
X_1 = \partial_x - \dfrac{y}{2} \partial_\tau, \ X_2 = \partial_y + \dfrac{x}{2} \partial_\tau.
\]
The Heisenberg group $\mathds{H}$ enjoys a structure of Lie group when equipped with the law
\[
(x, y, \tau) \cdot  (x', y', \tau') = (x + x', y + y', \tau + \tau' - \frac{1}{2} \mathrm{Im}[z \overline{z'}]),
\]
where $z := (x, y)$ and $z' := (x', y')$ are elements of $\mathds{R}^2$ that we will identify with $\mathds{C}$ for convenience ($\overline{\cdot}$ denotes the complex conjugation). The neutral element of this operation is $(0, 0, 0)$ and the inverse of $(x, y, \tau)$ is $(-x, -y, -\tau)$. There are many good references on the Heisenberg group, see \cite{heisenberg} or \cite{comprehensive2020} for example.

The Hamiltonian $H : \T^*(\He) \to \R$ is thus given by
\[
H(\lambda) = \dfrac{1}{2} \left(\left(v + \dfrac{\alpha x}{2}\right)^2 + \left(u - \dfrac{\alpha y}{2}\right)^2\right),
\]
for $\lambda = \left(q, u \mathrm{d}x|_q + v \mathrm{d}y|_q + \alpha \mathrm{d}\tau|_q \right)$ and $q = (x, y, \tau)$. In complex coordinates, we will write $z := x + i y$ and $w := u + i v$. Since the sub-Riemannian structure is left-invariant, we may choose $q = (0, 0, 0)$ for the rest of this section, without loss of generality.

We can explicitly solve Hamilton's equations to find the expression of a normal geodesic starting from` $q = (0, 0, 0)$ of $\He$ with initial covector $\lambda_0 = u_0 \diff x|_q + v_0 \diff y|_q + \alpha_0 \diff \tau|_q \in \T^*_q(\He)$:
\[
\lambda(t) = \left\{
\begin{array}{ll}
	z(t) & = \dfrac{1}{i \alpha_0} w_{0} (\me^{i \alpha_0 t} - 1) \\
	\tau(t) & = \dfrac{1}{2 \alpha_0^2} |w_{0}|^2 (\alpha_0 t - \sin(\alpha_0 t)) \\
	w(t) & = \dfrac{1}{2} w_{0} (\me^{i \alpha_0 t} + 1) \\
	\alpha(t) &= \alpha_0
\end{array}
\right.
\]

The symplectic moving frame $(E_1, E_2, E_3, F_1, F_2, F_3) = (\partial_u, \partial_v, \partial_\alpha, \partial_x, \partial_y, \partial_\tau)$ along $\lambda(t)$ induced by the global coordinates of $\mathds{H}^n$ satisfies \cref{movingframe}. The Jacobi fields along $\gamma$ are thus determined by the differential equation
\[
\begin{pmatrix}
	\dot{p} \\
	\dot{x}
\end{pmatrix}
=
\begin{pmatrix}
	- A(t)^T & R(t) \\
	B(t) & A(t)
\end{pmatrix}
\begin{pmatrix}
	p \\
	x
\end{pmatrix}.
\]

 A computation in the Heisenberg group with the chosen symplectic frame shows that the block matrices in the above differential equation are given by
\[
R(t) = 
\begin{pmatrix}
	- \dfrac{\alpha_0^2}{4} & 0 & 0 \\
	0 & - \dfrac{\alpha_0^2}{4} & 0 \\ 
	0 & 0 & 0 \\
\end{pmatrix},
\]
\[
A(t) = 
\begin{pmatrix}
	0 & - \dfrac{\alpha_0}{2} & 0 \\
	\dfrac{\alpha_0}{2} & 0 & 0 \\ 
	-\dfrac{v_0 - 3 (v_0 \cos(\alpha_0 t) + u_0 \sin(\alpha_0 t))}{4} & \dfrac{x_0 - 3 (x_0 \cos(\alpha_0 t) - v_0 \sin(\alpha_0 t))}{4} & 0 \\
\end{pmatrix}
\]
and the symmetric matrix
\[
B(t) = 
\begin{pmatrix}
	1 & 0 & -\dfrac{u_0 - u_0 \cos(\alpha_0 t) + v_0 \sin(\alpha_0 t)}{2 \alpha_0}
	) \\
	0 & 1 & -\dfrac{v_0 - v_0 \cos(\alpha_0 t) - u_0 \sin(\alpha_0 t)}{2 \alpha_0} \\ 
	\star & \star & \dfrac{4 |w_0|^2 (1 - \cos(\alpha_0 t))}{8 \alpha_0^2} \\
\end{pmatrix}.
\]
Solving this ordinary differential equation yields the general form of a Jacobi field $\mathcal{J}(t) = \sum_{i = 1}^{n} p_i(t) E_i(t) + x_i(t) F_i(t)$ along $\lambda(t)$ with initial condition $(p(0), x(0))$. After some calculations, we find
\begin{equation}
\label{jacobifieldsheis}
\begin{pmatrix}
	p(t) \\
	x(t)
\end{pmatrix}
=
M(t)
\begin{pmatrix}
	p_0 \\
	x_0
\end{pmatrix}
=
\begin{pmatrix}
	M_1(t) & M_2(t) \\
	M_3(t) & M_4(t)
\end{pmatrix}
\begin{pmatrix}
	p(0) \\
	x(0)
\end{pmatrix},
\end{equation}
where
\[
M(t) := 
\begin{pmatrix}
	\dfrac{1 + \cos(\alpha_0 t)}{2} & - \dfrac{\sin(\alpha_0 t)}{2} & f_1(t) & -\dfrac{\alpha_0 \sin(\alpha_0 t)}{4} & \dfrac{\alpha_0 (1 - \cos(\alpha_0 t))}{4} & 0 \\[2ex]
	\dfrac{\sin(\alpha_0 t)}{2} & \dfrac{1 + \cos(\alpha_0 t)}{2} & f_2(t) & \dfrac{\alpha_0 (1 - \cos(\alpha_0 t))}{4} & -\dfrac{\alpha_0 \sin(\alpha_0 t)}{4} & 0 \\[2ex]
	0 & 0 & 1 & 0 & 0 & 0 \\[2ex]
	\dfrac{\sin(\alpha_0 t)}{\alpha_0 } & -\dfrac{1 - \cos(\alpha_0 t)}{\alpha_0 } & f_3(t) & \dfrac{1 + \cos(\alpha_0 t)}{2} & - \dfrac{\sin(\alpha_0 t)}{2} & 0 \\[2ex]
	\dfrac{1 - \cos(\alpha_0 t)}{\alpha_0 } & \dfrac{\sin(\alpha_0 t)}{\alpha_0 } & f_4(t) & \dfrac{\sin(\alpha_0 t)}{2} & \dfrac{1 + \cos(\alpha_0 t)}{2} & 0 \\[2ex]
	f_5(t) & f_6(t) & f_7(t) & f_8(t) & f_9(t) & 1
\end{pmatrix}
\]
and
\[
\begin{cases}
	f_1(t) = \dfrac{- 2 t v_0 \cos(\alpha_0 t) - 2 t u_0 \sin(\alpha_0 t)}{4}; \\
	f_2(t) = \dfrac{2 t u_0 \cos(\alpha_0 t) - 2 t v_0 \sin(\alpha_0 t)}{4}; \\
	f_3(t) = \dfrac{- v_0 + (v_0 - t \alpha_0 u_0) \cos(t \alpha_0) + (u_0 + t \alpha_0 v_0) \sin(t \alpha_0)}{\alpha_0^2}; \\
	f_4(t) = \dfrac{- u_0 + (u_0 + t \alpha_0 v_0) \cos(\alpha_0 t) + (u_0 - t \alpha_0 u_0) \sin(\alpha_0 t)}{\alpha_0^2}; \\
	f_5(t) = \dfrac{2 t u_0 - 2 u_0 \sin(t \alpha_0)}{2 \alpha_0^2}; \\
	f_6(t) = \dfrac{2 t v_0 - 2 v_0 \sin(t \alpha_0)}{2 \alpha_0^2}; \\
	f_7(t) = |w_0|^2\dfrac{-\alpha_0 t(1 + \cos(\alpha_0 t))  + 2 \sin(\alpha_0 t)}{2 \alpha_0^3};\\
	f_8(t) = \dfrac{\alpha_0 t v_0 + u_0 (1 - \cos(\alpha_0 t))}{2 \alpha_0}; \\
	f_9(t) = \dfrac{- \alpha_0 t u_0 + v_0 (1 - \cos(\alpha_0 t))}{2 \alpha_0}.
\end{cases}
\] 

With this explicit description of sub-Riemannian Jacobi fields along normal geodesics of $\mathds{H}$, we are now able to study conjugates vectors. Alternatively, this could also be done by computing the determinant of the sub-Riemannian exponential map. We use the Jacobi fields characterisation of the kernel of $\exp_q$ to illustrate our work.

\begin{proposition}
	\label{conjlocusH}
	The covector $\lambda_0 = u_0 \diff x|_q + v_0 \diff y|_q + \alpha_0 \diff \tau|_q \in \T^*_q(\He)$ is a conjugate covector of $q = (0, 0, 0) \in \mathds{H}$ with $H(\lambda_0) \neq 0$ if and only if $\alpha_0 \sin(\alpha_0) + 2 \cos(\alpha_0) - 2 = 0$ and $\alpha_0 \neq 0$. Furthermore, the conjugate vectors are all of order one, they form a submanifold of $\T^*_q(\He)$ of dimension $2$ and
	\[
	\Kern(\diff_{\lambda_0} \exp_q) = 
	\left\{
	\begin{array}{ll}
		\vspan\left\{\begin{pmatrix}
			- v_0 \\
			u_0 \\
			0
		\end{pmatrix} \right\} & \text{ if } \sin(\alpha_0) = 0 \\
		\vspan\left\{\begin{pmatrix}
			-\frac{v_0 (\cos(\alpha_0) - 1)}{4} \\
			\frac{u_0 (\cos(\alpha_0) - 1)}{4} \\
			1
		\end{pmatrix} \right\} & \text{ if } \sin(\alpha_0) \neq 0
	\end{array}
	\right.
	\]
\end{proposition}

\begin{proof}
	We have seen in \cref{jacobifields} that
	\begin{align*}
		\Kern(\diff_{\lambda_0} \exp_q) & = \left\{ \sum_{k = 1}^3 p_k(1) E_k(1) \mid (p(t), x(t)) \text{ satisfies \eqref{jacobifieldsheis} with } x(0) = x(1) = 0 \right\} \\
		& \cong \left\{ \nabla J(1) \mid J \in \mathscr{J}_{0, 1}(\gamma_{q, \lambda_0}) \right\}.
	\end{align*}
	A covector $\lambda_0$ will be conjugate if the kernel of the $3 \times 3$ bottom left block matrix of $M(1)$ is non trivial. Furthermore, the image of those vectors under $M_1(1)$ will give $\Kern(\diff_{\lambda_0} \exp_q)$.
	
	If $\alpha_0$ tends to $0$, the block is similar to
	\[
	\begin{pmatrix}
		1 & 0 & - \frac{1}{2}v_0 \\
		0 & 1 & \frac{1}{2}u_0 \\
		0 & 0 & \frac{1}{12}|w_0|^2
	\end{pmatrix}.
	\]
	Since we assume that $H(\lambda_0) \neq 0$, the matrix above has a trivial kernel and $\alpha_0 = 0$ does not produce a conjugate vector.
	
	We assume now that $\alpha_0 \neq 0$. If $\sin(\alpha_0) = 0$ and $\cos(\alpha_0) = 1$, then we have
	\[
	\begin{pmatrix}
		0 & 0 & u_0 \\
		0 & 0 & v_0 \\
		u_0 & v_0 & \frac{-1 }{\alpha_0}|w_0|^2
	\end{pmatrix}.
	\]
	Since $H(\lambda_0) \neq 0$, the numbers $u_0$ and $v_0$ can not vanish at the same time. This case thus yields a conjugate vector of dimension 1 with $\Kern(\diff_{\lambda_0} \exp_q) = \vspan\{(- v_0, u_0, 0)\}$.
	
	If $\sin(\alpha_0) = 0$ and $\cos(\alpha_0) = - 1$, the block matrix $M_3(1)$ is similar to
	\[
	\begin{pmatrix}
		0 & 1 & \frac{u_0}{2} - \frac{v_0}{\alpha_0} \\
		1 & 0 & -\frac{v_0}{2} - \frac{u_0}{\alpha_0} \\
		0 & 0 & \frac{1}{2 \alpha_0^2} H(\lambda_0)
	\end{pmatrix}
	\]
	which has a trivial kernel. 
	
	If $\sin(\alpha_0) \neq 0$, the matrix is equivalent to
	\[
	\begin{pmatrix}
		1 & 0 & \star \\
		0 & 1 & \star \\
		0 & 0 & H(\lambda_0) \Bigl[\alpha_0 \cot\Bigl(\dfrac{\alpha_0}{2}\Bigr) - 2 \Bigr] \\
	\end{pmatrix}
	\]
	which has a non trivial kernel if only if $\alpha_0 \sin(\alpha_0) + 2 \cos(\alpha_0) - 2 = 0$. In this case, we obtain 
	\[
	\Kern(\diff_{\lambda_0} \exp_q) = \vspan\left\{
	\begin{pmatrix}
		-\frac{v_0 (\cos(\alpha_0) - 1)}{4} \\
		\frac{u_0 (\cos(\alpha_0) - 1)}{4} \\
		1
	\end{pmatrix}
	\right\}.
	\]
	Altogether, we finally observe that the collection of conjugate vectors consists of $(u_0, v_0, \alpha_0) \in \T^*_q(\He)$ such that $\alpha_0 \sin(\alpha_0) + 2 \cos(\alpha_0) - 2 = 0$, $\alpha_0 \neq 0$ and $(u_0, v_0) \neq (0, 0)$. They form planes in $\T^*_q(\He)$, parallel to the plane $\mathrm{span}\left\{\partial_u, \partial_v\right\}$, where the covector $\alpha_0 \partial_\alpha$ has been removed.
\end{proof}

\subsection{Local non injectivity of the exponential map of the Heisenberg group}

The set of all conjugate covectors to $q$ is called the conjugate locus at $q$ and we denote it by $\mathrm{Conj}_q(\mathds{H})$.

The structure of the conjugate locus being established, we can finally prove Morse-Littauer's theorem for the Heisenberg group by a novel application of Warner's approach from \cite[Theorem 3.4.]{warner1965} in sub-Riemannian geometry. We expect this method to work in other classes of sub-Riemannian manifolds.

\begin{proof}[Proof of \cref{thm:noninjH}.]
	Let $\lambda_0 \in \mathscr{A}_q \setminus H_q^{-1}(0)$ be a conjugate vector of $q \in \mathds{H}$. By \cref{conjlocusH}, we can thus choose a two-dimensional open connected submanifold $C$ of $\mathrm{Conj}_q(\mathds{H})$ in $\T^*_q(\mathds{H})$ containing $\lambda_0$. In particular, the conjugate vectors in $C$ have all order 1. We also write $C^0$ (resp. $C^1$) for the set of covectors $\overline{\lambda}_0 \in C$ such that 
	\[
	\dim \bigl[ \Kern(\diff_{\overline{\lambda}_0} \exp_q) \cap \T_{\overline{\lambda}_0}(\mathrm{Conj}_q(\mathds{H})) \bigr] = 0 \text{ (resp. = 1)}.
	\]
	It is easily seen that the set $C^1$ (resp. $C^0$) corresponds to the case $\sin(\alpha_0) = 0$ (resp. $\sin(\alpha_0) \neq 0$). Both $C^1$ and $C^0$ are thus open sets.
	
	\textbf{Case 1 : $\lambda_0 \in C^1$.}
	
	The subspaces $\Kern(\diff_{\overline{\lambda}_0} \exp_q)$ with $\overline{\lambda_0} \in C^1$ form a one-dimensional and involutive smooth distribution of $C^1$. This is because it corresponds to the distribution induced by the kernels of $\diff_{\overline{\lambda}_0} \exp_q$. By Frobenius theorem, there exists a unique integral manifold passing through $\lambda_0$. This is a one dimensional connected submanifold $N$ of $C^1$ such that $\T_{\overline{\lambda}_0}(N) = \Kern(\diff_{\overline{\lambda}_0} \exp_q)$ for all $\overline{\lambda}_0 \in N$. We then have that the restriction of $\exp_q$ to $N$ satisfies $\diff_{\overline{\lambda}_0} (\exp_q|_N) = 0$ for every $\overline{\lambda}_0 \in N$. Since $N$ is connected, this implies that the sub-Riemannian exponential map maps every elements of $N$ into a single point and hence $\exp_q$ is not injective in any neighbourhood of $\lambda_0 \in C^1$.
	
	\textbf{Case 2 : $\lambda_0 \in C^0$.}
	
	Firstly, the hypothesis that $\lambda_0 \in C^0$ means that $\Kern(\diff_{\overline{\lambda}_0} \exp_q) \not\subseteq \T_{\overline{\lambda}_0}(\mathrm{Conj}_q(\mathds{H}))$, i.e. $\lambda_0$ is a \emph{fold singularity} of $\exp_q$. We also have that
	\begin{equation}
	\label{ddtexpH}
	\dfrac{\diff}{\diff t}\Big|_{t = 1} \Bigl[ \diff_{t \lambda_0} \exp_q \Bigr] = \dfrac{H(\lambda_0)}{2 \alpha_0^2}\Bigl[ 2 - (2 + \alpha_0^2) \cos(\alpha_0) \Bigr].
	\end{equation}
	We know that $H(\lambda_0) \neq 0$ since $\lambda_0 \in \mathscr{A}_q \setminus H_q^{-1}(0)$. Furthermore, the covector $\lambda_0 \in C^0$ satisfies $\alpha_0 \sin(\alpha_0) + 2 \cos(\alpha_0) - 2 = 0$, $\alpha_0 \neq 0$ and $\sin(\alpha_0) \neq 0$. This means that
	\[
	\alpha_0 = \frac{2 - 2 \cos(\alpha_0)}{\sin(\alpha_0)} \ \text{ and } \ 2 - (2 + \alpha_0^2) \cos(\alpha_0) = \frac{16}{\sin^2(\alpha_0)} \sin^6(\alpha_0/2)
	\]
	which vanishes if and only if $\alpha_0 = 2 \pi n$ for some $n \in \mathds{Z}$. This never happens when $\lambda_0 \in C^0$. Therefore, the expression \eqref{ddtexpH} is not zero. This means that $\lambda_0$ is a \emph{good singularity} of $\exp_q$. By Whitney's singularity theory \cite{whitney1955} (see also Warner \cite[Theorem 3.3. c)]{warner1965}), we deduce that there exists coordinate systems $(\mathcal{U}, \xi)$ and $(\mathcal{V}, \eta)$ around $\lambda_0$ and $\exp_q(\lambda_0)$ respectively, such that
	\begin{enumerate}[label=\normalfont(\roman*)]
		\item $\eta^k \circ \exp_q = \xi^k$ for all $k = 1, 2$;
		\item $\eta^3 \circ \exp_q = \xi^3 \cdot \xi^3 $.
	\end{enumerate}
	This normal form of the sub-Riemannian exponential map implies that $\exp_q$ can not be injective in any neighbourhood of $\lambda_0 \in C^0$.
\end{proof}

We have described a cotangent version of Warner's method, via normal forms, to prove the non-injectivity of the Heisenberg exponential map in a neighbourhood of a singularity. We also point the reader to \cite[Figure 6]{lerariorizzi} for a related qualitative picture of the exponential map of the Heisenberg group.





\begingroup
\setstretch{0.9}
\setlength\bibitemsep{0pt}
\printbibliography

@article {sachkov2008,
    AUTHOR = {Sachkov, Yu. L.},
     TITLE = {Conjugate points in the {E}uler elastic problem},
   JOURNAL = {J. Dyn. Control Syst.},
  FJOURNAL = {Journal of Dynamical and Control Systems},
    VOLUME = {14},
      YEAR = {2008},
    NUMBER = {3},
     PAGES = {409--439},
      ISSN = {1079-2724},
   MRCLASS = {49Q10 (49J15 74G60 74G65 74K10 93B29 93C10)},
  MRNUMBER = {2425306},
MRREVIEWER = {Andrey V. Sarychev},
       DOI = {10.1007/s10883-008-9044-x},
       URL = {https://doi.org/10.1007/s10883-008-9044-x},
}

@article {arnold1997,
    AUTHOR = {Arnol'd, V. I.},
     TITLE = {On a characteristic class entering into conditions of quantization},
   JOURNAL = {Funkcional. Anal. i Prilo\v{z}en.},
  FJOURNAL = {Akademija Nauk SSSR. Funkcional'nyi Analiz i ego Prilo\v{z}enija},
    VOLUME = {1},
      YEAR = {1967},
     PAGES = {1--14},
      ISSN = {0374-1990},
   MRCLASS = {57.50 (35.00)},
  MRNUMBER = {0211415},
MRREVIEWER = {B. L. Reinhart},
}

@book {heisenberg,
    AUTHOR = {Capogna, Luca and Danielli, Donatella and Pauls, Scott D. and
              Tyson, Jeremy T.},
     TITLE = {An introduction to the {H}eisenberg group and the
              sub-{R}iemannian isoperimetric problem},
    SERIES = {Progress in Mathematics},
    VOLUME = {259},
 PUBLISHER = {Birkh\"{a}user Verlag, Basel},
      YEAR = {2007},
     PAGES = {xvi+223},
   MRCLASS = {53C17 (22E30 30C65 32T27 32V15 49Q15)},
  MRNUMBER = {2312336},
MRREVIEWER = {Piotr Haj\l asz},
}

@article {warner1965,
    AUTHOR = {Warner, Frank W.},
     TITLE = {The conjugate locus of a {R}iemannian manifold},
   JOURNAL = {Amer. J. Math.},
  FJOURNAL = {American Journal of Mathematics},
    VOLUME = {87},
      YEAR = {1965},
     PAGES = {575--604},
      ISSN = {0002-9327},
   MRCLASS = {53.72},
  MRNUMBER = {208534},
MRREVIEWER = {N. J. Hicks},
       DOI = {10.2307/2373064},
       URL = {https://doi.org/10.2307/2373064},
}

@article {rosquist1983,
	AUTHOR = {Rosquist, Kjell},
	TITLE = {On the structure of space-time caustics},
	JOURNAL = {Comm. Math. Phys.},
	FJOURNAL = {Communications in Mathematical Physics},
	VOLUME = {88},
	YEAR = {1983},
	NUMBER = {3},
	PAGES = {339--355},
	ISSN = {0010-3616},
	MRCLASS = {58E10 (53C50 83C99)},
	MRNUMBER = {701922},
	MRREVIEWER = {J. K. Beem},
	URL = {http://projecteuclid.org/euclid.cmp/1103922381},
}

@article {szeghy2008,
	AUTHOR = {Szeghy, D.},
	TITLE = {On the conjugate locus of pseudo-{R}iemannian manifolds},
	JOURNAL = {Indag. Math. (N.S.)},
	FJOURNAL = {Koninklijke Nederlandse Akademie van Wetenschappen.
	Indagationes Mathematicae. New Series},
	VOLUME = {19},
	YEAR = {2008},
	NUMBER = {3},
	PAGES = {465--480},
	ISSN = {0019-3577},
	MRCLASS = {53C50 (53C20)},
	MRNUMBER = {2513063},
	MRREVIEWER = {Francisco J. Palomo},
	DOI = {10.1016/S0019-3577(08)80013-4},
	URL = {https://doi.org/10.1016/S0019-3577(08)80013-4},
}

@book {comprehensive2020,
	AUTHOR = {Agrachev, Andrei and Barilari, Davide and Boscain, Ugo},
	title={A Comprehensive Introduction to Sub-Riemannian Geometry},
	SERIES = {Cambridge Studies in Advanced Mathematics},
	VOLUME = {181}, 
	PUBLISHER = {Cambridge University Press, Cambridge},
	YEAR = {2020},
	ISBN = {978-1-108-47635-5},
	MRCLASS = {53C17},
	MRNUMBER = {3971262},
	MRREVIEWER = {Luca Rizzi},
	DOI={10.1017/9781108677325}, 
}

@incollection {agrachev2008,
	AUTHOR = {Agrachev, A. A.},
	TITLE = {Geometry of optimal control problems and {H}amiltonian
	systems},
	BOOKTITLE = {Nonlinear and optimal control theory},
	SERIES = {Lecture Notes in Math.},
	VOLUME = {1932},
	PAGES = {1--59},
	PUBLISHER = {Springer, Berlin},
	YEAR = {2008},
	MRCLASS = {49N60 (37J05 37N35 49N90 58E25 93B29)},
	MRNUMBER = {2410710},
	MRREVIEWER = {Matthias Kawski},
	DOI = {10.1007/978-3-540-77653-6_1},
	URL = {https://doi.org/10.1007/978-3-540-77653-6_1},
}

@article {subriemjacobi,
	AUTHOR = {Barilari, Davide and Rizzi, Luca},
	TITLE = {On {J}acobi fields and a canonical connection in
	sub-{R}iemannian geometry},
	JOURNAL = {Arch. Math. (Brno)},
	FJOURNAL = {Universitatis Masarykianae Brunensis. Facultas Scientiarum
	Naturalium. Archivum Mathematicum},
	VOLUME = {53},
	YEAR = {2017},
	NUMBER = {2},
	PAGES = {77--92},
	ISSN = {0044-8753},
	MRCLASS = {53C17 (53B15 53B21)},
	MRNUMBER = {3672782},
	MRREVIEWER = {Pei Biao Zhao},
	DOI = {10.5817/AM2017-2-77},
	URL = {https://doi.org/10.5817/AM2017-2-77},
}

@book {piccione,
	AUTHOR = {Piccione, Paolo and Tausk, Daniel Victor},
	TITLE = {A student's guide to symplectic spaces, {G}rassmannians and
	{M}aslov index},
	SERIES = {Publica\c{c}\~{o}es Matem\'{a}ticas do IMPA. [IMPA Mathematical
	Publications]},
	PUBLISHER = {Instituto de Matem\'{a}tica Pura e Aplicada (IMPA), Rio de
	Janeiro},
	YEAR = {2008},
	PAGES = {xiv+301},
	ISBN = {978-85-244-0283-8},
	MRCLASS = {53D12 (53D05)},
	MRNUMBER = {2463212},
	MRREVIEWER = {Erasmo Caponio},
}

@article {morse,
	author = {Morse, Marston and Littauer, S. B.},
	title = {A Characterization of Fields in the Calculus of Variations},
	volume = {18},
	number = {12},
	pages = {724--730},
	year = {1932},
	doi = {10.1073/pnas.18.12.724},
	publisher = {National Academy of Sciences},
	issn = {0027-8424},
	URL = {https://www.pnas.org/content/18/12/724},
	eprint = {https://www.pnas.org/content/18/12/724.full.pdf},
	journal = {Proceedings of the National Academy of Sciences}
}

@article {whitney1955,
	AUTHOR = {Whitney, Hassler},
	TITLE = {On singularities of mappings of euclidean spaces. {I}.
	{M}appings of the plane into the plane},
	JOURNAL = {Ann. of Math. (2)},
	FJOURNAL = {Annals of Mathematics. Second Series},
	VOLUME = {62},
	YEAR = {1955},
	PAGES = {374--410},
	ISSN = {0003-486X},
	MRCLASS = {56.0X},
	MRNUMBER = {73980},
	MRREVIEWER = {R. Bott},
	DOI = {10.2307/1970070},
	URL = {https://doi.org/10.2307/1970070},
}

@article {lerariorizzi,
	AUTHOR = {Lerario, A. and Rizzi, L.},
	TITLE = {How many geodesics join two points on a contact
	sub-{R}iemannian manifold?},
	JOURNAL = {J. Symplectic Geom.},
	FJOURNAL = {The Journal of Symplectic Geometry},
	VOLUME = {15},
	YEAR = {2017},
	NUMBER = {1},
	PAGES = {247--305},
	ISSN = {1527-5256},
	MRCLASS = {53C17 (53C22)},
	MRNUMBER = {3652078},
	MRREVIEWER = {Enrico Le Donne},
	DOI = {10.4310/JSG.2017.v15.n1.a7},
	URL = {https://doi.org/10.4310/JSG.2017.v15.n1.a7},
}
\endgroup

\end{document}

Exponential map of the Heisenberg group for arbitrary $p$
//////////////////////////////

\section{Exponential map of the Heisenberg group} \label{examples}

In this section, we study the three dimensional Heisenberg group and prove \cref{thm:noninjH}. 

\subsection{Geometry of the Heisenberg group}

The Heisenberg group $\He$ is the sub-Riemannian structure, defined on $\R^3$, that is generated by the global vector fields
\[
X_1 = \partial_x - \dfrac{y}{2} \partial_\tau, \ X_2 = \partial_y + \dfrac{x}{2} \partial_\tau.
\]
The Heisenberg group $\mathds{H}$ enjoys a structure of Lie group when equipped with the law
\[
(x, y, \tau) \cdot  (x', y', \tau') = (x + x', y + y', \tau + \tau' - \frac{1}{2} \mathrm{Im}[z \overline{z'}]),
\]
where $z := (x, y)$ and $z' := (x', y')$ are elements of $\mathds{R}^2$ that we will identify with $\mathds{C}$ for convenience ($\overline{\cdot}$ denotes the complex conjugation). The neutral element of this operation is $(0, 0, 0)$ and the inverse of $(x, y, \tau)$ is $(-x, -y, -\tau)$. There are many good references on the Heisenberg group, see \cite{heisenberg} or \cite{comprehensive2020} for example.

The Hamiltonian $H : \T^*(\He) \to \R$ is thus given by
\[
H(\lambda) = \dfrac{1}{2} \left(\left(v + \dfrac{\alpha x}{2}\right)^2 + \left(u - \dfrac{\alpha y}{2}\right)^2\right),
\]
for $\lambda = \left(p, u \mathrm{d}x|_p + v \mathrm{d}y|_p + \alpha \mathrm{d}\tau|_p \right)$ and $p = (x, y, \tau)$.

We can explicitly solve Hamilton's equations to find the expression of a normal geodesic starting from an arbitrary point $p = (x_0, y_0, \tau_0)$ of $\He$ with initial covector $\lambda_0 = u_0 \diff x|_p + v_0 \diff y|_p + \alpha_0 \diff \tau|_p \in \T^*_p(\He)$:
\[
\lambda(t) = \left\{
\begin{array}{ll}
	z(t) & = \dfrac{1}{i \alpha_0} w_{0} (\me^{i \alpha_0 t} - 1) + \dfrac{1}{2} z_{0} (\me^{i \alpha_0 t} + 1) \\
	\tau(t) & = \tau_0 + \dfrac{\mathrm{Im}[\overline{z_{0}} w_{0}]}{2} t + \dfrac{1}{2 \alpha_0} \mathrm{Re}[\overline{z_{0}} w_{0}] (1 - \cos(\alpha_0 t)) \\
	& \qquad \qquad + \dfrac{1}{8} |z_{0}|^2 (\alpha_0 t + \sin(\alpha_0 t)) + \dfrac{1}{2 \alpha_0^2} |w_{0}|^2 (\alpha_0 t - \sin(\alpha_0 t)) \\
	w(t) & = \dfrac{1}{2} w_{0} (\me^{i \alpha_0 t} + 1) + \dfrac{\alpha_0}{4} i z_{0} (\me^{i \alpha_0 t} - 1) \\
	\alpha(t) &= \alpha_0
\end{array}
\right.
\]
To simplify the notations, we also set
\[
\xi_0 := u_0 - \dfrac{\alpha_0 y_0}{2}, \ \widetilde{\xi_0} := u_0 + \dfrac{\alpha_0 y_0}{2}, \ \eta_0 := v_0 + \dfrac{\alpha_0 x_0}{2} \text{ and } \widetilde{\eta_0} := v_0 - \dfrac{\alpha_0 x_0}{2}.
\]

The symplectic moving frame $(E_1, E_2, E_3, F_1, F_2, F_3) = (\partial_u, \partial_v, \partial_\alpha, \partial_x, \partial_y, \partial_\tau)$ along $\lambda(t)$ induced by the global coordinates of $\mathds{H}^n$ satisfies \cref{movingframe}. The Jacobi fields along $\gamma$ are thus determined by the differential equation
\[
\begin{pmatrix}
	\dot{p} \\
	\dot{x}
\end{pmatrix}
=
\begin{pmatrix}
	- A(t)^T & R(t) \\
	B(t) & A(t)
\end{pmatrix}
\begin{pmatrix}
	p \\
	x
\end{pmatrix}.
\]

After some computations, we find that in the Heisenberg group, and for the chosen symplectic frame, the block matrices in the above differential equation are given by
\[
R(t) = 
\begin{pmatrix}
	- \dfrac{\alpha_0^2}{4} & 0 & 0 \\
	0 & - \dfrac{\alpha_0^2}{4} & 0 \\ 
	0 & 0 & 0 \\
\end{pmatrix},
\]
\[
A(t) = 
\begin{pmatrix}
	0 & - \dfrac{\alpha_0}{2} & 0 \\
	\dfrac{\alpha_0}{2} & 0 & 0 \\ 
	-\dfrac{\widetilde{\eta_0} - 3 (\eta_0 \cos(\alpha_0 t) + \xi_0 \sin(\alpha_0 t))}{4} & \dfrac{\widetilde{\xi_0} - 3 (\xi_0 \cos(\alpha_0 t)) - \eta_0 \sin(\alpha_0 t)}{4} & 0 \\
\end{pmatrix}
\]
and the symmetric matrix
\[
B(t) = 
\begin{pmatrix}
	1 & 0 & -\dfrac{\widetilde{\xi_0} - \xi_0 \cos(\alpha_0 t) + \eta_0 \sin(\alpha_0 t)}{2 \alpha_0}
	) \\
	0 & 1 & -\dfrac{\widetilde{\eta_0} - \eta_0 \cos(\alpha_0 t) - \xi_0 \sin(\alpha_0 t)}{2 \alpha_0} \\ 
	\star & \star & \dfrac{4 |w_0|^2 (1 - \cos(\alpha_0 t)) + \alpha_0^2 |z_0|^2 (1 + \cos(\alpha_0 t)) + 4 \alpha_0 R[\overline{z_{0}} w_{0}] \sin(\alpha_0 t)}{8 \alpha_0^2} \\
\end{pmatrix}.
\]
Solving this ordinary differential equation yields the general form of a Jacobi field $\mathcal{J}(t) = \sum_{i = 1}^{n} p_i(t) E_i(t) + x_i(t) F_i(t)$ along $\lambda(t)$ with initial condition $(p(0), x(0))$. After some calculations, we find
\begin{equation}
\label{jacobifieldsheis}
\begin{pmatrix}
	p(t) \\
	x(t)
\end{pmatrix}
=
M(t)
\begin{pmatrix}
	p_0 \\
	x_0
\end{pmatrix}
=
\begin{pmatrix}
	M_1(t) & M_2(t) \\
	M_3(t) & M_4(t)
\end{pmatrix}
\begin{pmatrix}
	p(0) \\
	x(0)
\end{pmatrix},
\end{equation}
where
\[
M(t) := 
\begin{pmatrix}
	\dfrac{1 + \cos(\alpha_0 t)}{2} & - \dfrac{\sin(\alpha_0 t)}{2} & f_1(t) & -\dfrac{\alpha_0 \sin(\alpha_0 t)}{4} & \dfrac{\alpha_0 (1 - \cos(\alpha_0 t))}{4} & 0 \\[2ex]
	\dfrac{\sin(\alpha_0 t)}{2} & \dfrac{1 + \cos(\alpha_0 t)}{2} & f_2(t) & \dfrac{\alpha_0 (1 - \cos(\alpha_0 t))}{4} & -\dfrac{\alpha_0 \sin(\alpha_0 t)}{4} & 0 \\[2ex]
	0 & 0 & 1 & 0 & 0 & 0 \\[2ex]
	\dfrac{\sin(\alpha_0 t)}{\alpha_0 } & -\dfrac{1 - \cos(\alpha_0 t)}{\alpha_0 } & f_3(t) & \dfrac{1 + \cos(\alpha_0 t)}{2} & - \dfrac{\sin(\alpha_0 t)}{2} & 0 \\[2ex]
	\dfrac{1 - \cos(\alpha_0 t)}{\alpha_0 } & \dfrac{\sin(\alpha_0 t)}{\alpha_0 } & f_4(t) & \dfrac{\sin(\alpha_0 t)}{2} & \dfrac{1 + \cos(\alpha_0 t)}{2} & 0 \\[2ex]
	f_5(t) & f_6(t) & f_7(t) & f_8(t) & f_9(t) & 1
\end{pmatrix}
\]
and
\[
\begin{cases}
	f_1(t) = \dfrac{y_0 - (2 t \eta_0 + y_0) \cos(\alpha_0 t) - (2 t \xi_0 + x_0) \sin(\alpha_0 t)}{4}; \\
	f_2(t) = \dfrac{- x_0 + (2 t \xi_0 + x_0) \cos(\alpha_0 t) - (2 t \eta_0 + y_0) \sin(\alpha_0 t)}{4}; \\
	f_3(t) = \dfrac{- v_0 + (v_0 - t \alpha_0 \widetilde{\xi_0}) \cos(t \alpha_0) + (u_0 + t \alpha_0 \eta_0) \sin(t \alpha_0)}{\alpha_0^2}; \\
	f_4(t) = \dfrac{- u_0 + (u_0 + t \alpha_0 \eta_0) \cos(\alpha_0 t) + (u_0 - t \alpha_0 \widetilde{\xi_0}) \sin(\alpha_0 t)}{\alpha_0^2}; \\
	f_5(t) = \dfrac{2 t \xi_0 + \alpha_0 x_0 (1 - \cos(\alpha_0 t)) - 2 u_0 \sin(t \alpha_0)}{2 \alpha_0^2}; \\
	f_6(t) = \dfrac{2 t \eta_0 + \alpha_0 y_0 (1 - \cos(\alpha_0 t)) - 2 v_0 \sin(t \alpha_0)}{2 \alpha_0^2}; \\
	f_7(t) = \dfrac{1}{8 \alpha_0^3} \Bigl[\alpha_0 \bigl[t(\alpha_0^2 |z_0|^2 - 4 |w_0|^2)(1 + \cos(\alpha_0 t)) - 4 \ z_0 \cdot w_0 (1 - \cos(\alpha_0 t))\bigr]\\
	\qquad \qquad \qquad \qquad \qquad + 4 (2 |w_0|^2 + t \alpha_0^2 \ z_0 \cdot w_0 ) \sin(\alpha_0 t)\Bigr];\\
	f_8(t) = \dfrac{2 \alpha_0 t \eta_0 + 2 u_0 (1 - \cos(\alpha_0 t)) + \alpha_0 x_0 \sin(\alpha_0 t)}{4 \alpha_0}; \\
	f_9(t) = \dfrac{- 2 \alpha_0 t \xi_0 + 2 v_0 (1 - \cos(\alpha_0 t)) + \alpha_0 y_0 \sin(\alpha_0 t)}{4 \alpha_0}.
\end{cases}
\] 

With this explicit description of sub-Riemannian Jacobi fields along normal geodesics of $\mathds{H}$, we are now able to study conjugates vectors. Alternatively, this could also be done by computing the determinant of the sub-Riemannian exponential map. We use the Jacobi fields characterisation of the kernel of $\exp_p$ to illustrate our work.

\begin{proposition}
	\label{conjlocusH}
	The covector $\lambda_0 = u_0 \diff x|_p + v_0 \diff y|_p + \alpha_0 \diff \tau|_p \in \T^*_p(\He)$ is a conjugate covector of $p = (x_0, y_0, \tau_0) \in \mathds{H}$ with $H(\lambda_0) \neq 0$ if and only if $\alpha_0 \sin(\alpha_0) + 2 \cos(\alpha_0) - 2 = 0$ and $\alpha_0 \neq 0$. Furthermore, the conjugate vectors are all of order one, they form a submanifold of $\T^*_p(\He)$ of dimension $2$ and
	\[
	\Kern(\diff_{\lambda_0} \exp_p) = 
	\left\{
	\begin{array}{ll}
		\vspan\left\{\begin{pmatrix}
			- \eta_0 \\
			\xi_0 \\
			0
		\end{pmatrix} \right\} & \text{ if } \sin(\alpha_0) = 0 \\
		\vspan\left\{\begin{pmatrix}
			-\frac{\eta_0 (\cos(\alpha_0) - 1) - 2 y_0}{4} \\
			\frac{\xi_0 (\cos(\alpha_0) - 1) - 2 x_0}{4} \\
			1
		\end{pmatrix} \right\} & \text{ if } \sin(\alpha_0) \neq 0
	\end{array}
	\right.
	\]
\end{proposition}

\begin{proof}
	We have seen in \cref{jacobifields} that
	\begin{align*}
		\Kern(\diff_{\lambda_0} \exp_p) & = \left\{ \sum_{k = 1}^3 p_k(1) E_k(1) \mid (p(t), x(t)) \text{ satisfies \eqref{jacobifieldsheis} with } x(0) = x(1) = 0 \right\} \\
		& \cong \left\{ \nabla J(1) \mid J \in \mathscr{J}_{0, 1}(\gamma_{p, \lambda_0}) \right\}.
	\end{align*}
	A covector $\lambda_0$ will be conjugate if the kernel of the $3 \times 3$ bottom left block matrix of $M(1)$ is non trivial. Furthermore, the image of those vectors under $M_1(1)$ will give $\Kern(\diff_{\lambda_0} \exp_p)$.
	
	If $\alpha_0$ tends to $0$, the block is similar to
	\[
	\begin{pmatrix}
		1 & 0 & - \frac{1}{2}(v_0 + y_0) \\
		0 & 1 & \frac{1}{2}(u_0 + x_0) \\
		0 & 0 & \frac{1}{12}|w_0|^2
	\end{pmatrix}.
	\]
	Since we assume that $H(\lambda_0) \neq 0$, the matrix above has a trivial kernel and $\alpha_0 = 0$ does not produce a conjugate vector.
	
	We assume now that $\alpha_0 \neq 0$. If $\sin(\alpha_0) = 0$ and $\cos(\alpha_0) = 1$, then we have
	\[
	\begin{pmatrix}
		0 & 0 & \xi_0 \\
		0 & 0 & \eta_0 \\
		\xi_0 & \eta_0 & \frac{-4 |w_0|^2 + \alpha_0^2 |z_0|^2}{4 \alpha_0}
	\end{pmatrix}.
	\]
	Since $H(\lambda_0) \neq 0$, the numbers $\xi_0$ and $\eta_0$ can not vanish at the same time. This case thus yields a conjugate vector of dimension 1 with $\Kern(\diff_{\lambda_0} \exp_p) = \vspan\{(- \eta_0, \xi_0, 0)\}$.
	
	If $\sin(\alpha_0) = 0$ and $\cos(\alpha_0) = - 1$, the block matrix $M_3(1)$ is similar to
	\[
	\begin{pmatrix}
		0 & 1 & \frac{\xi_0}{2} - \frac{v_0}{\alpha_0} \\
		1 & 0 & -\frac{\eta_0}{2} - \frac{u_0}{\alpha_0} \\
		0 & 0 & \frac{1}{2 \alpha_0^2} H(\lambda_0)
	\end{pmatrix}
	\]
	which has a trivial kernel. 
	
	If $\sin(w_0) \neq 0$, the matrix is equivalent to
	\[
	\begin{pmatrix}
		1 & 0 & \star \\
		0 & 1 & \star \\
		0 & 0 & H(\lambda_0) \Bigl[\alpha_0 \cot\Bigl(\dfrac{\alpha_0}{2}\Bigr) - 2 \Bigr] \\
	\end{pmatrix}
	\]
	which has a non trivial kernel if only if $\alpha_0 \sin(\alpha_0) + 2 \cos(\alpha_0) - 2 = 0$. In this case, we obtain 
	\[
	\Kern(\diff_{\lambda_0} \exp_p) = \vspan\left\{
	\begin{pmatrix}
		-\frac{\eta_0 (\cos(\alpha_0) - 1) - 2 y_0}{4} \\
		\frac{\xi_0 (\cos(\alpha_0) - 1) - 2 x_0}{4} \\
		1
	\end{pmatrix}
	\right\}.
	\]
	We finally observe that the collection of conjugate vectors consists of $(u_0, v_0, \alpha_0) \in \T^*_p(\He)$ such that $\sin(\alpha_0)\left(\alpha_0 \cot(\alpha_0/2) - 2 \right) = 0$, $\alpha_0 \neq 0$ and $(u_0, v_0) \neq \frac{\alpha_0}{2}(y_0, - x_0)$. They form planes in $\T^*_p(\He)$, parallel to the plane $\mathrm{span}\left\{\partial_u, \partial_v\right\}$, where the covector $\frac{\alpha_0}{2}(y_0 \partial_u - x_0 \partial_v) + \alpha_0 \partial_\alpha$ has been removed.
\end{proof}

\subsection{Local non injectivity of the exponential map of the Heisenberg group}

The structure of the conjugate locus being established, we can finally prove Morse-Littauer's theorem for the Heisenberg group by a novel application of Warner's approach from \cite[Theorem 3.4.]{warner1965} in sub-Riemannian geometry. We expect this method to work in other classes of sub-Riemannian manifolds.

\begin{proof}[Proof of \cref{thm:noninjH}.]
	Let $\lambda_0 \in \mathscr{A}_p \setminus H_p^{-1}(0)$ be a conjugate vector of $p \in \mathds{H}$. We denote by $\mathrm{Conj}_p(\mathds{H})$ the conjugate locus at $p$. By \cref{conjlocusH}, we can thus choose a two-dimensional open connected submanifold $C$ of $\mathrm{Conj}_p(\mathds{H})$ in $\T^*_p(\mathds{H})$ containing $\lambda_0$. In particular, the conjugate vectors in $C$ have all order 1. We also write $C^0$ (resp. $C^1$) for the set of covectors $\overline{\lambda}_0 \in C$ such that 
	\[
	\dim \bigl[ \Kern(\diff_{\overline{\lambda}_0} \exp_p) \cap \T_{\overline{\lambda}_0}(\mathrm{Conj}_p(\mathds{H})) \bigr] = 0 \text{ (resp. = 1)}.
	\]
	It is easily seen that the set $C^1$ (resp. $C^0$) corresponds to the case $\sin(\alpha_0) = 0$ (resp. $\sin(\alpha_0) \neq 0$). Both $C^1$ and $C^0$ are thus open sets.
	
	\textbf{Case 1 : $\lambda_0 \in C^1$.}
	
	The subspaces $\Kern(\diff_{\overline{\lambda}_0} \exp_p)$ with $\overline{\lambda_0} \in C^1$ form a one-dimensional and involutive smooth distribution of $C^1$. This is because it corresponds to the distribution induced by the kernels of $\diff_{\overline{\lambda}_0} \exp_p$. By Frobenius theorem, there exists a unique integral manifold passing through $\lambda_0$. This is a one dimensional connected submanifold $N$ of $C^1$ such that $\T_{\overline{\lambda}_0}(N) = \Kern(\diff_{\overline{\lambda}_0} \exp_p)$ for all $\overline{\lambda}_0 \in N$. We then have that the restriction of $\exp_p$ to $N$ satisfies $\diff_{\overline{\lambda}_0} (\exp_p|_N) = 0$ for every $\overline{\lambda}_0 \in N$. Since $N$ is connected, this implies that the sub-Riemannian exponential map maps every elements of $N$ into a single point and hence $\exp_p$ is not injective in any neighbourhood of $\lambda_0 \in C^1$.
	
	\textbf{Case 2 : $\lambda_0 \in C^0$.}
	
	Firstly, the hypothesis that $\lambda_0 \in C^0$ means that $\Kern(\diff_{\overline{\lambda}_0} \exp_p) \not\subseteq \T_{\overline{\lambda}_0}(\mathrm{Conj}_p(\mathds{H}))$, i.e. $\lambda_0$ is a \emph{fold singularity} of $\exp_p$. We also have that
	\begin{equation}
	\label{ddtexpH}
	\dfrac{\diff}{\diff t}\Big|_{t = 1} \Bigl[ \diff_{t \lambda_0} \exp_p \Bigr] = \dfrac{H(\lambda_0)}{2 \alpha_0^2}\Bigl[ 2 - (2 + \alpha_0^2) \cos(\alpha_0) \Bigr].
	\end{equation}
	We know that $H(\lambda_0) \neq 0$ since $\lambda_0 \in \mathscr{A}_p \setminus H_p^{-1}(0)$. Furthermore, the covector $\lambda_0 \in C^0$ satisfies $\alpha_0 \sin(\alpha_0) + 2 \cos(\alpha_0) - 2 = 0$, $\alpha_0 \neq 0$ and $\sin(\alpha_0) \neq 0$. Therefore, the expression \eqref{ddtexpH} is not zero. This means that $\lambda_0$ is a \emph{good singularity} of $\exp_p$. By Whitney's singularity theory \cite{whitney1955} (see also Warner \cite[Theorem 3.3. c)]{warner1965}), we deduce that there exists coordinate systems $(\mathcal{U}, \xi)$ and $(\mathcal{V}, \eta)$ around $\lambda_0$ and $\exp_p(\lambda_0)$ respectively, such that
	\begin{enumerate}[label=\normalfont(\roman*)]
		\item $\eta^k \circ \exp_p = \xi^k$ for all $k = 1, 2$;
		\item $\eta^3 \circ \exp_p = \xi^3 \cdot \xi^3 $.
	\end{enumerate}
	This normal form of the sub-Riemannian exponential map implies that $\exp_p$ can not be injective in any neighbourhood of $\lambda_0 \in C^0$.
\end{proof}

We have described a cotangent version of Warner's method, via normal forms, to prove the non-injectivity of the Heisenberg exponential map in a neighbourhood of a singularity. We also point the reader to \cite[Figure 6]{lerariorizzi} for a related qualitative picture of the exponential map of the Heisenberg group.

\section*{Data Availability Statement}

Data sharing not applicable to this article as no datasets were generated or analysed during the current study.